\documentclass[11pt]{amsart}
\usepackage{amsmath,amssymb,amsthm,geometry,hyperref,tikz}
\usetikzlibrary{arrows.meta}
\definecolor{strandred}{HTML}{D53A42}
\definecolor{strandblue}{HTML}{2469BC}
\definecolor{strandgreen}{HTML}{218A50}
\definecolor{strandturquoise}{HTML}{079AAB}
\definecolor{twistpurple}{HTML}{854AB2}
\hypersetup{colorlinks=true,citecolor=blue,linkcolor=blue,urlcolor=blue}
\newtheorem{theorem}{Theorem}[section]
\newtheorem{proposition}[theorem]{Proposition}
\newtheorem{lemma}[theorem]{Lemma}
\newtheorem{corollary}[theorem]{Corollary}
\theoremstyle{remark}
\newtheorem{example}[theorem]{Example}
\newtheorem{remark}[theorem]{Remark}

\newcommand{\Z}{\mathbb Z}
\newcommand{\Q}{\mathbb Q}
\newcommand{\Int}{\operatorname{Int}}
\newcommand{\Mod}{\operatorname{Mod}}
\newcommand{\Conf}{\operatorname{Conf}}

\title{Fibered-knot braids and symplectic fake quadrics}

\author{Anar Akhmedov}
\address{School of Mathematics, University of Minnesota, Minneapolis, MN 55455, USA}
\email{akhmedov@umn.edu}

\author{B. Doug Park}
\address{Department of Pure Mathematics, University of Waterloo, Waterloo, ON, N2L 3G1, Canada}
\email{bdpark@uwaterloo.ca}

\author{Sai-Kee Yeung}
\address{Department of Mathematics, Purdue University, West Lafayette, IN 47907, USA}
\email{yeungs@purdue.edu}

\date{October 8, 2026}



\begin{document}

\begin{abstract}
We construct aspherical symplectic four-manifolds with the rational
cohomology ring of $S^2\times S^2$ by summing ruled surfaces along
braided multisections. 
For each $m\geq1$, we form two symplectic sums with the same gluing map in the
chosen surface markings and normal framings, using a stabilized braid
for the trefoil in one sum and for the torus knot $T(2,7)$ in the other. 
Both sums have integral first 
homology $(\Z/2)^2\oplus\Z/(4m)$. Their fundamental groups admit
respectively $36$ and $24$ epimorphisms onto $S_3$, so the sums are
nonhomeomorphic. We also vary the gluing while keeping the
multisections fixed. In each of five product, twisted, and hybrid
cases, this gives an infinite sequence in which the first homology
groups have strictly increasing finite orders. All but finitely many
members of each sequence admit no complex structure.
\end{abstract}

\maketitle

\section{Introduction}

Akhmedov and Yeung~\cite{AY} constructed symplectic analogues of
fake quadrics by summing ruled surfaces along cyclic multisections.
Their construction gives four-manifolds with the rational cohomology
ring of $S^2\times S^2$. Here we ask whether the closed sum can
distinguish two braided multisections that define the same covering
of the base. We show that changing the braid can change the
fundamental group of the sum even when the surface gluing is fixed
and the integral first homology groups agree.

In $\Sigma_2\times S^2$, Smith~\cite[Proposition~1.2]{Smith} used
hyperelliptic monodromy and branched covers to construct infinitely
many nonisotopic genus-five symplectic surfaces in the class
$4[\Sigma_2\times\{\mathrm{pt}\}]$. Hays and Park~\cite{HaysPark}
extended this construction to even genus. We obtain representatives
of the same class from the four-strand braids
\[
  \beta_k=\sigma_1^{2k+1}\sigma_2\sigma_3,\qquad k\geq1.
\]
Each $\beta_k$ is obtained by two positive Markov stabilizations of
the standard two-strand braid for the torus knot $T(2,2k+1)$. Following the
construction in Section~\ref{sec:braids}, we let the four points in
a sphere fiber move by $\beta_k$ around the generator $a_1$ of
$\pi_1(\Sigma_2)$ and remain fixed around the other standard
generators. This defines a connected degree-four multisection
$F_{\beta_k}\subset\Sigma_2\times S^2$. Since every $\beta_k$
induces the same four-cycle, the coverings
$F_{\beta_k}\to\Sigma_2$ are isomorphic over $\Sigma_2$.
For the other summand we use the degree-two multisection
$F_{3,2}\subset\Sigma_3\times S^2$ defined by the half-twist
$\sigma_1\in B_2$. Both summing surfaces have genus five and
self-intersection zero.

We use the surface generators and normal framings of
Section~\ref{sec:fixed-gluing}. For each $m\geq1$, the based gluing
$\psi_m:F_{\beta_k}\to F_{3,2}$ sends the marked generators of the
source surface group to the same words in the target generators
for every $k$. In the specified framings, the boundary gluing is
$(x,z)\mapsto(\psi_m(x),\overline z)$ for $z\in S^1$. Denote the
resulting product sum by
\[
  Z_{k,m}=Z_{k,\psi_m}:=
  (\Sigma_2\times S^2)
  \mathbin{\#}_{F_{\beta_k}=F_{3,2},\,\psi_m}
  (\Sigma_3\times S^2).
\]
The choices $k=1$ and $k=3$ correspond to the right-handed trefoil
and $T(2,7)$, respectively.

\begin{theorem}\label{thm:main-fixed-gluing}
For every $m\geq1$, the manifolds $Z_{1,m}$ and $Z_{3,m}$ are
closed, aspherical, minimal symplectic four-manifolds with the
rational cohomology ring of $S^2\times S^2$. For $k\in\{1,3\}$,
\[
  e(Z_{k,m})=4,\qquad \sigma(Z_{k,m})=0,\qquad
  H_1(Z_{k,m};\Z)\cong(\Z/2)^2\oplus\Z/(4m).
\]
For each fixed $m$, the groups $\pi_1(Z_{1,m})$ and
$\pi_1(Z_{3,m})$ admit respectively $36$ and $24$ epimorphisms
onto $S_3$. In particular, $Z_{1,m}$ and $Z_{3,m}$ are
nonhomeomorphic. For all sufficiently large $m$, neither manifold
admits a complex structure.
\end{theorem}

The distinction comes from the action of $\beta_k$ on the
fundamental group of the four-punctured sphere. On first homology
this action only permutes the meridian classes. To compute the
fundamental group of the sum, we determine how the boundary
inclusions map a meridian and the pushed-off surface generators
into each exterior group. This gives the presentation of
$\pi_1(Z_{k,m})$ in Proposition~\ref{prop:explicit-closed-group}.
For every $k$ and $m$, six epimorphisms onto $S_3$ kill the
meridian and eighteen send it to a three-cycle. The number that
send it to a transposition is twelve when $3\nmid k$ and zero
when $3\mid k$ (Proposition~\ref{prop:s3-distinction}).

We next consider what happens when the two multisections are fixed
and the surface gluing varies. Section~\ref{sec:five-patterns}
treats three pairs of degrees in product sphere bundles, namely
$(2,3)$, $(2,4)$, and $(3,3)$, together with two further $(2,4)$
pairs. The twisted pair uses two nontrivial sphere bundles, and
the hybrid pair uses one product and one nontrivial bundle.
In each case we obtain a sequence of aspherical minimal symplectic
four-manifolds with the rational cohomology ring of $S^2\times S^2$.
Their integral first homology groups have strictly increasing finite
orders, so the members of each sequence are pairwise nonhomeomorphic
(Theorem~\ref{thm:five-families}). If one of these manifolds admits
a complex structure, it is a minimal surface of general type with
$K^2=8$ and $\chi(\mathcal O)=1$. There are only finitely many
differentiable types of such surfaces. The increasing homology
orders therefore rule out complex structures for all sufficiently
large $m$. In each of the three $(2,4)$ cases, this threshold is
independent of the braid parameter $k$
(Theorem~\ref{thm:noncomplex}).

The covering maps of cyclic multisections were studied in
\cite{AkhmedovMultisections}, and the intersection lattices of
higher mixed sums in \cite{AkhmedovHigherMixed}. For the three
product pairs, replacing a multisection by its braided representative
preserves the covering map up to isomorphism and the order of the
meridian in the first homology of the exterior. With adapted normal
framings, the integral first homology formulas of
\cite[Theorem~4.5]{AY} therefore apply. In the twisted and hybrid
pairs, the degree-four covering has dihedral permutation monodromy
of order eight. We show in Section~\ref{sec:five-braids} that its
pushforward on integral first homology has image of index two.
The cyclic degree-four pushforward used in
\cite[Theorem~8.2 and Proposition~9.4]{AY} has image of index four.
Proposition~\ref{prop:unbounded-gluing} applies directly to the
covering maps in both settings and gives gluings for which the
orders of the first homology groups are unbounded. This argument
requires no fiberedness assumption on a braid closure.

Earlier applications of fibered knots to symplectic sums include
\cite{AkhmedovCohomology,AkhmedovSmall,AkhmedovParkSmall,AkhmedovParkOdd,AkhmedovParkSakalli}.
For knots of genus one or two, a second way to obtain braids uses
the hyperelliptic quotient of the capped knot fiber
\cite{AkhmedovBraidedCover}. The monodromy descends to a sphere
with marked branch points, and each orbit of length $p$ gives a
degree-$p$ component of the branch locus. For the trefoil and the
figure-eight knot, the orbit lengths are $1$ and $3$
(Section~\ref{sec:four-strand}). For the even-indexed Stallings
twists of the square knot in the family of Park and Yun~\cite{ParkYun},
there are two orbits of length three (Section~\ref{sec:genus-two}).
For Misev's genus-two knots~\cite{Misev}, there is a fixed point
and an orbit of length five. In the final section we use the
resulting degree-five multisections to sum two copies of
$\Sigma_2\times S^2$ along genus-six surfaces. For every knot in
the family, the gluing can be chosen so that the sum has $e=12$,
$\sigma=0$, and $b_1=0$ (Propositions~\ref{prop:five}
and~\ref{prop:twelve}).

\section{Braids and complements}
\label{sec:braids}

We use the commutator convention $[u,v]=uvu^{-1}v^{-1}$ and we compose braid actions by
$\Phi_{\alpha\beta}=\Phi_\alpha\circ\Phi_\beta$. Our convention for
bundle monodromy is that conjugation by a lifted base generator acts by
the corresponding automorphism of the fiber group. All peripheral words
below use these conventions. 
Throughout, $\Sigma_g$ denotes a closed
connected oriented surface of genus $g$, and $B_p$ is the Artin braid
group on $p$ strands in a disk. The group $\Mod(S^2,p)$ consists of
orientation-preserving mapping classes preserving a set of $p$ marked
points; the individual points may be permuted. Hyperelliptic quotients
give elements of this group. 

Fix integers $g\geq1$, $p\geq1$, and fix a smoothly embedded closed
disk $D\subset S^2$, equipped with the orientation induced
from $S^2$.
Write $\Conf_p(\Int D)$ for the unordered configuration space of $p$
distinct points in $\Int D$, and choose a smooth loop representing
$\beta\in B_p$. We denote this loop also by $\beta$. 

Choose a basepoint $y_0\in \Sigma_g$ and 
a smooth map $u:\Sigma_g\to S^1$ whose induced homomorphism
\[
n=u_*:\pi_1(\Sigma_g,y_0)\longrightarrow\pi_1(S^1,1)\cong \Z
\]
satisfies $n(a_1)=1$, $n(b_1)=0$, and
$n(a_i)=n(b_i)=0$ for $i\geq2$, where
$a_1,b_1,\ldots,a_g,b_g$ are the standard generators of $\pi_1(\Sigma_g)$. 
Such a map $u$\/ exists: map the $a_1$-edge of the standard CW
decomposition of $\Sigma_g$ once around $S^1$ and all other
edges constantly. The attaching loop of the two-cell has
degree zero, so this map extends over $\Sigma_g$ and can
then be replaced by a smooth map in its homotopy class.

Let us define
\[
F_\beta=\{(s,z)\in\Sigma_g\times\Int D\mid z\in\beta(u(s))\}
\subset \Sigma_g\times S^2,
\]
and let $q_\beta:F_\beta\to\Sigma_g$ and
$r_\beta:F_\beta\to\Int D$ be the two projections. Thus $F_\beta$
is the pullback of the tautological $p$-sheeted cover of the configuration
space. 

\begin{proposition}\label{prop:braid-surface}
If the permutation of $\beta$ is a $p$-cycle, then $F_\beta$ is a
connected embedded surface and $q_\beta$ is an unbranched covering of
degree $p$. With the orientation induced by this covering,
\[
[F_\beta]=p[\Sigma_g\times\{\mathrm{pt}\}]\in 
H_2(\Sigma_g\times S^2;\mathbb Z),
\]
$[F_\beta]^2=0$, and $g(F_\beta)=p(g-1)+1$.
Moreover, if $\omega_g$ and $\omega_{S^2}$
are positive area forms, then $F_\beta$ is symplectic for every product
form
\[
\omega_\lambda=\lambda\,\operatorname{pr}_1^*\omega_g
                 +\operatorname{pr}_2^*\omega_{S^2} \quad(\lambda>0)
\]
on $\Sigma_g\times S^2$.
\end{proposition}

\begin{proof}
Locally the configuration consists of $p$ smooth disjoint branches, so
$F_\beta$ is embedded and $q_\beta$ is a local diffeomorphism. The
monodromy of this covering sends $a_1$ to the strand permutation and
all the other standard generators to the identity. 
Since the strand permutation of $\beta$ is a $p$-cycle, the
monodromy action on the sheets is transitive. Hence $F_\beta$
is connected, and
\[
2-2g(F_\beta)=p(2-2g).
\]
The projection $\Sigma_g\times D\to\Sigma_g$ is a homotopy equivalence
and sends $[F_\beta]$ to $p[\Sigma_g]$. This proves the homology formula
in $\Sigma_g\times D$, and hence in $\Sigma_g\times S^2$.

Because $dq_\beta$ is an isomorphism at each point, the vertical tangent
planes form a complement to $TF_\beta$ in the ambient tangent bundle.
Consequently the normal bundle is isomorphic to $r_\beta^*TD$, which is
trivial. In particular, $[F_\beta]^2=0$.

On each local branch, the sphere coordinate is a function of the single
parameter $u(s)$. Hence $dr_\beta$ has rank at most one and
$r_\beta^*\omega_{S^2}=0$. It follows that
\[
\omega_\lambda|_{F_\beta}=\lambda q_\beta^*\omega_g,
\]
which is positive for every $\lambda>0$.
\end{proof}

Choose a closed tubular neighborhood $\nu F_\beta$ consisting of small
vertical disks, and set
\begin{equation}\label{eq:M_beta}
M_\beta=(\Sigma_g\times S^2)\setminus\Int\nu F_\beta,
\qquad
P_p=S^2\setminus\coprod_{j=1}^p\Int D_j.
\end{equation}
We call $M_\beta$ the \emph{exterior}, or the \emph{complement}, of $F_\beta$. 
Its projection to $\Sigma_g$ is a bundle with fiber $P_p$. 
Fix a basepoint $z_0\in S^2\setminus D$. The constant section
$\Sigma_g\times\{z_0\}$ is disjoint from $\nu F_\beta$ and
therefore defines a section of $M_\beta\to\Sigma_g$.

Choose positively oriented meridians $x_1,\ldots,x_p$, based at $z_0$ by a standard system of paths, so that
\[
\pi_1(P_p,z_0)=\langle x_1,\ldots,x_p\mid x_1\cdots x_p=1\rangle.
\]
The braid motion fixes the chosen basepoint $z_0$ throughout, so the monodromy of
$M_\beta\to\Sigma_g$ acts by automorphisms of $\pi_1(P_p)$, rather
than only outer automorphisms. 
The monodromy of $q_\beta$ records the permutation of the
$p$ sheets, while the monodromy of $M_\beta\to\Sigma_g$
records the induced automorphism of $\pi_1(P_p,z_0)$.
Braids with the same strand permutation can induce different
automorphisms of this group.

Since $\pi_2(\Sigma_g)=0$, the homotopy sequence and the $z_0$ section give
a split exact sequence
\[
1\longrightarrow\pi_1(P_p,z_0)\longrightarrow\pi_1(M_\beta,(y_0,z_0))
\longrightarrow\pi_1(\Sigma_g,y_0)\longrightarrow1.
\]
The bundle is pulled back from the configuration loop along $u$, so
its based action factors through $n$. Thus we have 
\[
\pi_1(M_\beta)\cong
\pi_1(P_p)\rtimes_{\rho_\beta}\pi_1(\Sigma_g),
\]
where the based monodromy representation is
\[
\rho_\beta:\pi_1(\Sigma_g)\longrightarrow
\operatorname{Aut}(\pi_1(P_p,z_0)),\qquad
\rho_\beta(h)=\Phi_\beta^{\,n(h)}
\quad(h\in\pi_1(\Sigma_g)).
\]
In particular, only $a_1$ can act nontrivially among the standard base
generators. On first homology the action only permutes the meridians.

\begin{lemma}\label{lem:h1}
If the permutation of $\beta$ is a $p$-cycle, then
\[
H_1(M_\beta;\mathbb Z)
\cong H_1(\Sigma_g;\mathbb Z)\oplus\mathbb Z/p,
\]
where the $\mathbb Z/p$ summand is generated by $\mu$,
the homology class of a positively oriented meridian of
$F_\beta\subset\Sigma_g\times S^2$.
\end{lemma}

\begin{proof}
By the preceding split extension,
\[
\pi_1(M_\beta)\cong
\pi_1(P_p)\rtimes_{\rho_\beta}\pi_1(\Sigma_g).
\]
Abelianizing gives the direct sum of $H_1(\Sigma_g;\mathbb Z)$
and the quotient of $H_1(P_p;\mathbb Z)$ obtained by identifying
each class with its monodromy images. Now
\[
H_1(P_p;\mathbb Z)
=
\left\langle [x_1],\ldots,[x_p]\ \middle|\
[x_1]+\cdots+[x_p]=0\right\rangle.
\]
Since the strand permutation is a $p$-cycle, the monodromy
relations identify all $[x_j]$ with a single class $\mu$.
The resulting quotient therefore has the presentation
\[
\langle\mu\mid p\mu=0\rangle\cong\mathbb Z/p.
\]
The splitting of the base contribution is induced by the
constant $z_0$ section, and $\mu$ is the homology class of a
positively oriented meridian of $F_\beta$.
\end{proof}

We next compute the fundamental group of the complement $M_\beta$.
Here $\beta\in B_p$ is the braid used to construct $F_\beta$ from
the configuration map $\beta\circ u$. Recall that $u_*(a_1)=1$
and that $u_*$ vanishes on the other standard generators of
$\pi_1(\Sigma_g)$. The following presentation is valid for any
$\beta\in B_p$, whether or not its strand permutation is a
$p$-cycle.

\begin{proposition}\label{prop:group-presentation}
Let $g\geq1$ and $p\geq1$, and let $M_\beta$ be the complement
associated to $\beta\in B_p$ defined in\/ {\rm (\ref{eq:M_beta})}.
Then $G_{\beta,g}:=\pi_1(M_\beta,(y_0,z_0))$ has the presentation
\[
\begin{split}
G_{\beta,g}=\langle a_i,b_i,x_1,\ldots,x_p\mid\;&
x_1\cdots x_p=1,\quad
\prod_{i=1}^{g}[a_i,b_i]=1,\\
&a_1x_ja_1^{-1}=\Phi_\beta(x_j),\quad [b_1,x_j]=1,\\
&[a_i,x_j]=[b_i,x_j]=1\quad(2\leq i\leq g)\rangle.
\end{split}
\]
Here $a_i,b_i$ range over $1\leq i\leq g$, and all relations
involving $x_j$ are imposed for $1\leq j\leq p$. The generators
$a_i,b_i$ are lifted through the constant section $\Sigma_g \times \{z_0\}$, and the
$x_j$ are based meridians in the fiber $\{y_0\}\times P_p$.

The automorphism $\Phi_\beta$ is the based braid action on
$\pi_1(P_p,z_0)$. With the convention
$\Phi_{\alpha\beta}=\Phi_\alpha\circ\Phi_\beta$, it is computed
from the Artin rules
\[
\Phi_{\sigma_i}(x_i)=x_ix_{i+1}x_i^{-1},\qquad
\Phi_{\sigma_i}(x_{i+1})=x_i,\qquad
\Phi_{\sigma_i}(x_j)=x_j
\mathrm{ \ \ if \ } j\notin\{i,i+1\},
\]
using inverse automorphisms for negative braid generators.
\end{proposition}

\begin{proof}
Choose a reference fiber $\{y_0\}\times P_p$ and a standard system of based meridians
$x_1,\ldots,x_p$. These give
\[
\pi_1(P_p,z_0)=
\langle x_1,\ldots,x_p\mid x_1\cdots x_p=1\rangle.
\]
The displayed Artin rules describe the action of a positive
half twist on the meridians of the punctured disk. They preserve
the product $x_1\cdots x_p$, since
\[
\Phi_{\sigma_i}(x_i)\Phi_{\sigma_i}(x_{i+1})
=(x_ix_{i+1}x_i^{-1})x_i=x_ix_{i+1},
\]
and all other generators are fixed. Hence the Artin action
descends to $\pi_1(P_p)$, giving the automorphism $\Phi_\beta$.

By the split extension established above,
\[
\pi_1(M_\beta)\cong
\pi_1(P_p)\rtimes_{\rho_\beta}\pi_1(\Sigma_g).
\]
Lift the standard generators of $\pi_1(\Sigma_g,y_0)$ through the constant section $\Sigma_g\times\{z_0\}$
and denote their lifts by the same letters. Because the section
induces a homomorphism on fundamental groups, these lifts satisfy
\[
\prod_{i=1}^g[a_i,b_i]=1
\]
in $\pi_1(M_\beta)$.

The monodromy factors through $u_*$. Thus conjugation by $a_1$
acts on the fiber group by $\Phi_\beta$, while conjugation by
$b_1,a_i,b_i$ for $i\geq2$ acts trivially. These are precisely
the remaining displayed relations. A presentation of a
semidirect product consists of presentations of the two factors
together with the relations specifying this conjugation action,
so the displayed relations give a complete presentation of
$\pi_1(M_\beta)$.
\end{proof}

\begin{example}\label{ex:beta_k}
Consider $\beta=\beta_k=\sigma_1^{2k+1}\sigma_2\sigma_3$, with
$k\geq1$. The four-strand complement in $\Sigma_2\times S^2$ has
generators $x_1,\ldots,x_4,a,b,c,d$ and relations
\begin{equation}\label{eq:4-strand-relations}
\begin{aligned}
x_1x_2x_3x_4&=1,\qquad [a,b][c,d]=1,\\
ax_ja^{-1}&=\Phi_{\beta_k}(x_j),\qquad
[b,x_j]=[c,x_j]=[d,x_j]=1\quad(1\leq j\leq4).
\end{aligned}
\end{equation}
To make the Artin substitution explicit, put
\[
\Delta=x_1x_2,\qquad
y=\Delta^k(x_1x_2x_1^{-1})\Delta^{-k},\qquad
z=\Delta^kx_1\Delta^{-k}.
\]
The square $\Phi_{\sigma_1}^2$ conjugates both $x_1$ and $x_2$ by
$\Delta$, fixes $x_3,x_4$, and fixes $\Delta$. Consequently
$\Phi_{\sigma_1^{2k+1}}(x_1)=y$ and
$\Phi_{\sigma_1^{2k+1}}(x_2)=z$. Applying the factors in the prescribed
order gives
\begin{equation}\label{eq:Artin}
\begin{aligned}
\Phi_{\beta_k}(x_1)&=y,&
\Phi_{\beta_k}(x_2)&=zx_3z^{-1},\\
\Phi_{\beta_k}(x_3)&=zx_4z^{-1},&
\Phi_{\beta_k}(x_4)&=z.
\end{aligned}
\end{equation}

For the two-strand complement in $\Sigma_3\times S^2$, write $t=x_1$.
The sphere relation gives $x_2=t^{-1}$, and the Artin rule for a strand
exchange sends $t$ to $t^{-1}$. The presentation is therefore
\begin{equation}\label{eq:2-strand-group}
\left\langle t,a_i,b_i\ (1\leq i\leq3)\ \middle|\
\begin{gathered}
\prod_{i=1}^{3}[a_i,b_i]=1,\qquad a_1ta_1^{-1}=t^{-1},\\
[b_1,t]=1,\qquad [a_i,t]=[b_i,t]=1\quad(i=2,3)\;
\end{gathered}\right\rangle.
\end{equation}
The preceding four-strand action has permutation $(1\,2\,3\,4)$.
Abelianizing the two presentations gives, respectively,
$\Z^4\oplus\Z/4$ and $\Z^6\oplus\Z/2$, in agreement with
Lemma~\ref{lem:h1}.
\end{example}

We next determine the image of the boundary inclusion
$\partial\nu F_\beta\hookrightarrow M_\beta$. Assume that
$p\geq2$ and that the strand permutation of $\beta$ is a
$p$-cycle. Let $\widetilde y_0$ be the unique point of
\[
F_\beta\cap
\bigl(\{y_0\}\times\operatorname{Int}D_1\bigr),
\]
and let
\[
\chi_p:\pi_1(\Sigma_g,y_0)\longrightarrow\mathbb Z/p
\]
be the reduction of $n=u_*$ modulo $p$. A loop based at $y_0$
lifts to a loop based at $\widetilde y_0$ precisely when its
image under $n$ is divisible by $p$. Consequently, the covering
injection identifies
\[
(q_\beta)_*\pi_1(F_\beta,\widetilde y_0)=\ker\chi_p.
\]

We specify the basepoints and connecting path used for the
boundary homomorphism. Let $\delta$ be the path in the
reference fiber $\{y_0\}\times P_p$ from $(y_0,z_0)$ to
$\{y_0\}\times\partial D_1$ used to define the based meridian
$x_1$. Thus $x_1$ is represented by following $\delta$,
traversing $\{y_0\}\times\partial D_1$ once in the positive
normal direction, and returning along $\delta^{-1}$.

Since the normal bundle of $F_\beta$ is trivial, choose an
oriented normal framing whose boundary identification
\[
F_\beta\times S^1\cong\partial\nu F_\beta
\]
sends $(\widetilde y_0,1)$ to $\delta(1)$. Let
$\iota:F_\beta\times S^1\to M_\beta$ denote this identification
followed by the boundary inclusion. The inclusion, followed
by change of basepoint along $\delta$, gives the homomorphism
\[
\begin{aligned}
j_\beta:\pi_1\bigl(F_\beta\times S^1,(\widetilde y_0,1)\bigr)
&\longrightarrow
\pi_1\bigl(M_\beta,(y_0,z_0)\bigr),\\
[\gamma]
&\longmapsto
\bigl[\delta\cdot(\iota\circ\gamma)\cdot\delta^{-1}\bigr],
\end{aligned}
\]
where concatenated paths are traversed from left to right.
We call its image the \emph{peripheral subgroup}\/ of $F_\beta$.
By construction, the positive normal-circle generator,
represented by
\[
t\longmapsto(\widetilde y_0,e^{2\pi it}),
\qquad 0\leq t\leq1,
\]
maps to the chosen meridian $x_1$.

\begin{proposition}\label{prop:peripheral-presentation}
The boundary inclusion induces an injection
\[
j_\beta: \pi_1(F_\beta\times S^1)\longrightarrow G_{\beta,g}
\]
whose image is the centralizer $C_{G_{\beta,g}}(x_1)$. Projection to
$\pi_1(\Sigma_g)$ gives an exact sequence
\begin{equation}\label{eq:peripheral-exact-sequence}
1\longrightarrow\langle x_1\rangle
\longrightarrow C_{G_{\beta,g}}(x_1)
\longrightarrow\ker\chi_p\longrightarrow1.
\end{equation}
The chosen normal framing splits this sequence. In particular, the
peripheral subgroup is isomorphic to $\pi_1(F_\beta)\times\Z$.
\end{proposition}

\begin{proof}
Let $N=\pi_1(P_p)$ and let $J$ be the image of $j_\beta$.
The circle factor maps to $x_1$, and projection to the base restricts
on the other factor to the covering injection
$\pi_1(F_\beta)\hookrightarrow\pi_1(\Sigma_g)$.
If an element of the boundary group $\pi_1(F_\beta\times S^1)$ maps to the identity, this covering
injection shows that its surface component is trivial. It is therefore
a power of the normal circle. The group $N$ is free of rank $p-1$,
with $x_1$ as one of its free generators, and embeds in $G_{\beta,g}$.
Thus $x_1$ has infinite order, and that power is also trivial. This
proves injectivity. 

The boundary projection to $\Sigma_g$ factors as
$F_\beta\times S^1\to F_\beta\xrightarrow{q_\beta}\Sigma_g$.
Since $(q_\beta)_*$ is injective with image $\ker\chi_p$,
the projection of $J$ has image $\ker\chi_p$ and kernel
equal to the image of the normal-circle factor, namely
$\langle x_1\rangle$.
Since the circle factor is central
in the boundary group, $J\subset C_{G_{\beta,g}}(x_1)$.

The centralizer of a free generator in a free group is the cyclic
subgroup it generates. Hence
\[
C_{G_{\beta,g}}(x_1)\cap N=C_N(x_1)=\langle x_1\rangle.
\]
This includes $p=2$ case, when $N=\langle x_1\rangle$ is itself cyclic.

Let $v\in G_{\beta,g}$ project to $h\in\pi_1(\Sigma_g,y_0)$,
and let $\widehat h$ be the lift of $h$ induced by the constant
section. Then $v=w\widehat h$ for some $w\in N$, and hence
\[
vx_1v^{-1}=w\,\rho_\beta(h)(x_1)\,w^{-1}.
\]
If the permutation monodromy of $h$ sends the label $1$ to $k$,
then $\rho_\beta(h)(x_1)$ is conjugate in $N$ to $x_k$:
monodromy carries a positive meridian around the first puncture
to a positive meridian around the $k$th puncture, while the
choice of basing path accounts for the conjugation.
Thus $vx_1v^{-1}$ is also conjugate in $N$ to $x_k$.
If $v$ centralizes $x_1$, it follows that $x_1$ and $x_k$
are conjugate in $N$ and therefore represent the same class
in $H_1(P_p;\mathbb Z)$. Since the classes of the oriented
meridians $x_1,\ldots,x_p$ are distinct, we must have $k=1$.
Consequently, $h$ fixes the chosen sheet, which is equivalent
to $h\in\ker\chi_p$.

Now take $v\in C_{G_{\beta,g}}(x_1)$. Choose $j\in J$ with the
same base projection. Then $vj^{-1}\in N$ and centralizes $x_1$, so
$vj^{-1}\in\langle x_1\rangle\subset J$. Thus $v\in J$, proving that 
$J=C_{G_{\beta,g}}(x_1)$. Finally, the product structure supplied by the
normal framing gives the claimed splitting.
\end{proof}

We use the term ``normal connected sum'' for the smooth manifold
obtained by gluing the exteriors along their boundary circle
bundles, reversing the circle orientations. 
When this construction is performed with compatible symplectic data,
we call it a ``symplectic sum''. The fundamental-group and
asphericity arguments below apply to the underlying normal
connected sum and do not require a symplectic structure.

For the chosen normal framing of $F_\beta$, define
\[
\ell^{(\beta)}:\pi_1(F_\beta)\longrightarrow G_{\beta,g},
\qquad
\ell^{(\beta)}(\alpha)=j_\beta(\alpha,0).
\]
Thus
\[
j_\beta(\alpha,m)=\ell^{(\beta)}(\alpha)x_1^m.
\]
The map $\ell^{(\beta)}$ represents loops pushed into the
boundary with constant normal-circle coordinate. 
The composition of $\ell^{(\beta)}$ with the natural
projection $G_{\beta,g}\to\pi_1(\Sigma_g)$ equals
$(q_\beta)_*$. Thus, after identifying $\pi_1(F_\beta)$
with $\ker\chi_p$ via $(q_\beta)_*$, the homomorphism
$\ell^{(\beta)}$ gives a splitting of the exact sequence
(\ref{eq:peripheral-exact-sequence}) in Proposition~\ref{prop:peripheral-presentation}.

We want to express the lifts of elements of $\ker\chi_p$ to $C_{G_{\beta,g}}(x_1)$ in terms of the complement generators.  
Let $h=(q_\beta)_*(\alpha)\in \ker\chi_p$ and let $\widehat h$ be its lift
to the section $\Sigma_g\times\{z_0\}$. 
Since $n(h)$ is divisible by $p$,
the permutation monodromy of $h$ fixes the point
$\widetilde y_0$. Hence the induced automorphism preserves
the conjugacy class of $x_1$, so there exists
$U_h\in\pi_1(P_p,z_0)$ such that
\[
\widehat h\,x_1\,\widehat h^{-1}=\Phi_\beta^{\,n(h)}(x_1)=U_hx_1U_h^{-1}.
\]
Then $U_h^{-1}\widehat h$ centralizes $x_1$ and projects to $h$.
Any two such centralizing lifts differ by a power of $x_1$. 

By Proposition~\ref{prop:peripheral-presentation}, we have 
\[
\ell^{(\beta)}(\alpha)
=x_1^{m_h}U_h^{-1}\widehat h
\]
for a unique integer $m_h$. For the chosen $U_h$, this
integer records the meridian correction prescribed by
the framing.

\begin{example}
Recall that the four-strand braid $\beta_k$ in Example~\ref{ex:beta_k} 
induces the cycle $(1\,2\,3\,4)$ on the strand labels, so $a^4\in\ker\chi_4$. Its action is the
fourth iterate $\Phi_{\beta_k}^{4}$; an identity
\[
\Phi_{\beta_k}^{4}(x_1)=Ux_1U^{-1}
\]
gives the peripheral lift $U^{-1}a^4$, up to the meridian
power determined by the framing.
\end{example}

Now let $F_\beta\subset\Sigma_g\times S^2$ and
$F_\gamma\subset\Sigma_{g'}\times S^2$ be connected
multisections of degrees $p,q\geq2$ with common genus
\[
G=p(g-1)+1=q(g'-1)+1.
\]
On the second exterior, make the analogous choices of
basepoints, normal framing, and connecting path.
Write $\widetilde y_0\hspace{-3pt}{}'$ for the basepoint of $F_\gamma$,
$x'_1$ for its based meridian, and
\[
\ell^{(\gamma)}:
\pi_1(F_\gamma,\widetilde y_0\hspace{-3pt}{}')\longrightarrow G_{\gamma,g'},
\qquad
\ell^{(\gamma)}(\alpha)=j_\gamma(\alpha,0).
\]
Fix a basepoint $s_0\in\Sigma_G$ and choose
orientation-preserving diffeomorphisms
\[
(\Sigma_G,s_0)\longrightarrow(F_\beta,\widetilde y_0),
\qquad
(\Sigma_G,s_0)\longrightarrow(F_\gamma,\widetilde y_0\hspace{-3pt}{}'\hspace{1pt}).
\]
Using these identifications, we can regard $\ell^{(\beta)}$
and $\ell^{(\gamma)}$ as homomorphisms from
\[
\pi_1(\Sigma_G,s_0)=
\left\langle z_1,\ldots,z_{2G}\ \middle|\
\prod_{i=1}^{G}[z_{2i-1},z_{2i}]=1\right\rangle,
\]
and write $\ell_j^{(\beta)}=\ell^{(\beta)}(z_j)$
and $\ell_j^{(\gamma)}=\ell^{(\gamma)}(z_j)$.

Next we choose an orientation-preserving diffeomorphism
$\psi:\Sigma_G\to\Sigma_G$ and a smooth map
$\theta:\Sigma_G\to S^1$ satisfying
\[
\psi(s_0)=s_0,\qquad \theta(s_0)=1.
\]
In the boundary product coordinates supplied by the
framings and surface identifications, glue by
\[
(s,\zeta)\longmapsto
\bigl(\psi(s),\theta(s)\overline{\zeta}\bigr).
\]
This map reverses the normal-circle orientation and
sends the boundary basepoint $(s_0,1)$ on the first
side to $(s_0,1)$ on the second.

Let $f=\theta_*:H_1(\Sigma_G;\mathbb Z)\to\mathbb Z$
and $f_j=f([z_j])$.  
Seifert-Van Kampen theorem gives the
fundamental group of the resulting normal connected
sum as the quotient of $G_{\beta,g}*G_{\gamma,g'}$ by
\[
x_1=(x'_1)^{-1},
\qquad
\ell_j^{(\beta)}
=\ell^{(\gamma)}\bigl(\psi_*(z_j)\bigr)(x'_1)^{f_j}
\quad (1\leq j\leq2G).
\]
The first relation identifies the oppositely oriented
meridians; the remaining relations record the surface
gluing and the winding of $\theta$ along each generator.
Here $\ell^{(\gamma)}(\psi_*(z_j))$ is obtained by
substituting the framed lifts $\ell_i^{(\gamma)}$
into the word $\psi_*(z_j)$.

\begin{remark}
We can examine how much information the complement group
retains beyond the strand permutation. Two braids with the
same permutation differ by a pure braid, so we first
consider the action of pure braids on the fiber group.

For $p=2$, every pure braid acts trivially on
$\pi_1(P_2)\cong\mathbb Z$. For $p=3$, every pure braid
acts by an inner automorphism, since the pure mapping
class group of the three-punctured sphere is trivial;
see \cite[Section~2.2]{FarbMargalit}.

Consequently, if two three-strand braids $\beta,\gamma$
have the same permutation, then
\[
\Phi_\gamma=\operatorname{Inn}(w)\circ\Phi_\beta
\]
for some $w\in\pi_1(P_3)$. In the presentation of
$G_{\gamma,g}$, replacing $a_1$ by
$\widetilde a_1=w^{-1}a_1$ changes its action to
$\Phi_\beta$. It also preserves the surface relation:
because $b_1$ and all $a_i,b_i$ with $i\geq2$ commute
with $w$,
\[
[\widetilde a_1,b_1]\prod_{i=2}^{g}[a_i,b_i]
=
w^{-1}\left(\prod_{i=1}^{g}[a_i,b_i]\right)w
=1.
\]
This invertible substitution gives
$G_{\beta,g}\cong G_{\gamma,g}$, but need not preserve
the prescribed framed boundary maps. With four punctures,
pure braid monodromy can have a nontrivial outer action,
so the complement group can retain information beyond
the strand permutation.
\end{remark}

\begin{remark}
The trace calculation in the next section distinguishes
marked, fiber-preserving monodromies; it does not prove
ambient non-isotopy of the multisections. An ambient
isotopy identifies the resulting sums when the normal
framings and gluing maps are transported through it.
The infinite families constructed later vary the gluing
and require no assertion of ambient non-isotopy.
\end{remark}

Asphericity follows uniformly from the exterior bundles.
We state the result for arbitrary oriented sphere bundles
to include the twisted and hybrid sums considered below.

\begin{proposition}\label{prop:aspherical-sums}
Let $E_i\to\Sigma_{g_i}$, $i=1,2$, be oriented sphere
bundles with $g_i\geq1$, and let $F_i\subset E_i$ be
connected unbranched multisections of degrees $p_i\geq2$.
Their exteriors
\[
M_i=E_i\setminus\operatorname{Int}\nu F_i
\]
are aspherical, and the boundary inclusions induce
injections on fundamental groups. Every normal connected
sum along $F_1,F_2$ is aspherical. In particular, every
such symplectic sum is minimal.
\end{proposition}

\begin{proof}
The exterior $M_i$ is a bundle over $\Sigma_{g_i}$ with
fiber $P_{p_i}$. Both base and fiber are aspherical,
so the homotopy exact sequence gives
$\pi_n(M_i)=0$ for $n\geq2$. The boundary $\partial M_i$
is a circle bundle over $F_i$ and is likewise aspherical.

An element in the kernel of
$\pi_1(\partial M_i)\to\pi_1(M_i)$ projects trivially
to $\pi_1(\Sigma_{g_i})$, and hence to $\pi_1(F_i)$,
because $F_i\to\Sigma_{g_i}$ is a covering.
It is therefore a power of the meridian.
The fiber group $\pi_1(P_{p_i})$ injects into
$\pi_1(M_i)$, and the meridian has infinite order in
this fiber group. Thus the boundary homomorphism is
injective.  
Seifert-Van Kampen theorem identifies the sum group with the
amalgam of the exterior groups along the boundary group.

Let $Z$ denote the normal connected sum and $\widetilde Z$ its universal
cover. 
By the normal-form theorem for amalgamated products
\cite[Chapter~I, \S1.2, Remark after Theorem~1]{SerreTrees},
the exterior groups and the boundary group inject
into $\pi_1(Z)$.
Consequently, each component of the inverse image of an exterior or of the gluing
boundary is its universal cover and is therefore
contractible.

Form a graph whose vertices correspond to the lifted
exteriors and whose edges correspond to the lifted
boundary collars joining them. 
This incidence graph is the tree associated to the
amalgamated product
\cite[Chapter~I, \S4.1, Theorem~7]{SerreTrees}.
The union associated to any finite subtree is contractible: starting with
one lifted exterior, each additional edge attaches a
contractible lifted exterior along a contractible
lifted boundary. These attachments preserve
contractibility because the boundary inclusions
have collars.

The lifted decomposition is locally finite, since
each point has a neighborhood meeting at most two
lifted exteriors. Thus the image of any map
$S^n\to\widetilde Z$ meets only finitely many lifted
pieces and lies in the union associated to a finite
subtree. That union is contractible, so the map is
null-homotopic. Hence all homotopy groups of
$\widetilde Z$ vanish, and Whitehead's theorem implies
that $\widetilde Z$ is contractible.

Finally, a sphere of self-intersection $-1$ cannot be
null-homotopic. An aspherical symplectic four-manifold
therefore contains no exceptional sphere and is minimal.
\end{proof}

\section{Four-strand knot braids and sums with Euler characteristic four}
\label{sec:four-strand}

We describe two ways in which knots supply braid monodromies
for the constructions that will follow.  
The first starts with the monodromy of a genus-one fibered
knot and takes the elliptic quotient of its capped fiber. 
This will produce a section and a degree-three multisections in product sphere bundles
$\Sigma_g\times S^2$.  
In the symplectic sums of
Section~\ref{sec:five-braids}, these are paired either with
degree-two multisections or with other degree-three multisections,
with the base genera chosen so that the summing surfaces have
the same genus.

The second starts with a disk braid whose closure is a knot and uses
the braid directly; see Proposition~\ref{prop:braid-surface}. 
We use this second construction to obtain the degree-four multisections appearing in Proposition~\ref{prop:four-sum} and in the fixed-gluing calculation of Section~\ref{sec:fixed-gluing}.

Throughout this section, $3_1$ denotes the right-handed trefoil
$T(2,3)$, and $4_1$ denotes the figure-eight knot. Their fibers are
genus-one surfaces with one boundary component. Choose each knot
monodromy to be the identity near the boundary. Capping the boundary
circle with a disk and extending the monodromy by the identity gives
a diffeomorphism of a closed torus. Identify each capped fiber with
$T=\mathbb R^2/\mathbb Z^2$. With respect to suitable symplectic bases
of $H_1(T;\mathbb Z)$, the capped monodromies are represented by
\[
M_{3_1}=\begin{pmatrix}1&1\\-1&0\end{pmatrix},
\qquad
M_{4_1}=\begin{pmatrix}0&-1\\1&3\end{pmatrix};
\]
see \cite[proof of Proposition~5.14]{BurdeZieschang}.
Their traces are $1$ and $3$, respectively.

The corresponding linear torus automorphisms commute with the
elliptic involution
\[
\iota:T\longrightarrow T,\qquad \iota(x)=-x.
\]
The quotient $T/\langle\iota\rangle$ is a sphere, and the quotient map
is a double branched covering. Its four ramification points form the
two-torsion subgroup
\[
T[2]=\bigl(\tfrac12\mathbb Z^2\bigr)/\mathbb Z^2.
\]
Under the identification
\[
H_1(T;\mathbb Z/2)\cong(\mathbb Z/2)^2
\xrightarrow{\;\cong\;}T[2],
\qquad [v]\longmapsto[v/2],
\]
the actions on these points are given by the reductions
\[
\overline M_{3_1}=\begin{pmatrix}1&1\\1&0\end{pmatrix},
\qquad
\overline M_{4_1}=\begin{pmatrix}0&1\\1&1\end{pmatrix}.
\]
Each fixes the origin and cyclically permutes the three nonzero
elements. Thus each induced sphere monodromy fixes one branch point
and cyclically permutes the other three.  

For $K\in\{3_1,4_1\}$, the mapping torus of the capped monodromy is
diffeomorphic to $S^3_0(K)$, the three-manifold obtained by zero-framed
surgery on $K$. Indeed, the boundary of a knot fiber represents the
preferred longitude, so the meridional disks of the attached solid
torus cap the boundary circles of the fibers. 

To obtain a smooth branched covering of mapping tori, choose
coordinates near the ramification points and their images in which
the covering is given by $w=z^2$. Within the induced mapping class
of the four-pointed sphere, choose a diffeomorphism that is a rotation
from a small coordinate disk about each branch point to a disk about
its image. Such a representative can be obtained by an isotopy
preserving the marked set, since only the centers of the disks are
marked. Lift this isotopy starting at the linear torus automorphism.
The resulting lift is isotopic to the capped knot monodromy and is
smooth at every ramification point: in the chosen coordinates, a
rotation $w\mapsto e^{i\theta}w$ lifts locally to
$z\mapsto\pm e^{i\theta/2}z$. The quotient map therefore induces a
smooth double branched covering between the corresponding mapping tori.

Since the sphere monodromies are orientation preserving, they are
isotopic to the identity after forgetting the branch points, so their
mapping tori are diffeomorphic to $S^1\times S^2$. Taking the product
with an additional circle therefore gives double branched coverings
\[
S^1\times S^3_0(K)\longrightarrow T^2\times S^2.
\]

In each case the branch locus is $B_1\sqcup B_3$. The component $B_1$
is the product of the additional $S^1$ factor with the closed curve
traced by the branch point corresponding to the origin in the torus $T$. 
The component $B_3$ is the product of this same $S^1$ factor
with the single closed curve traced by the three branch points
corresponding to the nonzero elements of $T[2]$. These curves project
to the base circle with degrees one and three, respectively. Thus
$B_1$ is a section and $B_3$ is a connected degree-three multisection
of $T^2\times S^2\to T^2$, and both components are tori.

The fixed section also provides a convenient product identification.
The sphere monodromy is isotopic to the identity while fixing the
branch point corresponding to the origin, provided the other three
points are allowed to move. Choose the product identification using
such an isotopy. Then $B_1$ is a constant section, and the other three
points trace a loop of configurations in its complementary plane.
This loop is contained in a disk and hence defines a three-strand
disk braid $\gamma$. Its permutation is a three-cycle. In particular,
$B_3$ lies in $T^2\times D$ for a disk $D\subset S^2$, and
\[
[B_3]=3[T^2\times\{\mathrm{pt}\}],\qquad B_3^2=0.
\]
Pulling the configuration loop back along the map
$u:\Sigma_g\to S^1$ of Section~\ref{sec:braids},
Proposition~\ref{prop:braid-surface} gives a connected symplectic
degree-three multisection $F_\gamma\subset\Sigma_g\times S^2$ with
\[
[F_\gamma]=3[\Sigma_g\times\{\mathrm{pt}\}],\qquad
F_\gamma^2=0,\qquad
g(F_\gamma)=3(g-1)+1.
\]
These are the degree-three inputs to be used later.  

We now turn to the degree-four construction. 
Consider the following four-strand braids:  
\[
\beta_{3_1}=\sigma_1^3\sigma_2\sigma_3,
\qquad
\beta_{4_1}=(\sigma_1\sigma_2^{-1})^2\sigma_3,
\qquad
\beta_{5_1}=\sigma_1^5\sigma_2\sigma_3.
\]
These are obtained from 
$\sigma_1^3\in B_2$, $(\sigma_1\sigma_2^{-1})^2\in B_3$,
and $\sigma_1^5\in B_2$ by two, one, and two positive Markov
stabilizations, respectively. Each stabilization adds one
strand and one positive crossing while preserving the isotopy
type of the closure. Thus their closures are the right-handed
trefoil, the figure-eight knot, and $T(2,5)$, respectively.

Their closures are connected, so their strand permutations are
four-cycles. In particular, all four points belong to one orbit in
this construction, whereas the preceding elliptic quotients have a
fixed point and a three-point orbit. These are different sphere
monodromies, even when the source knot is the same.

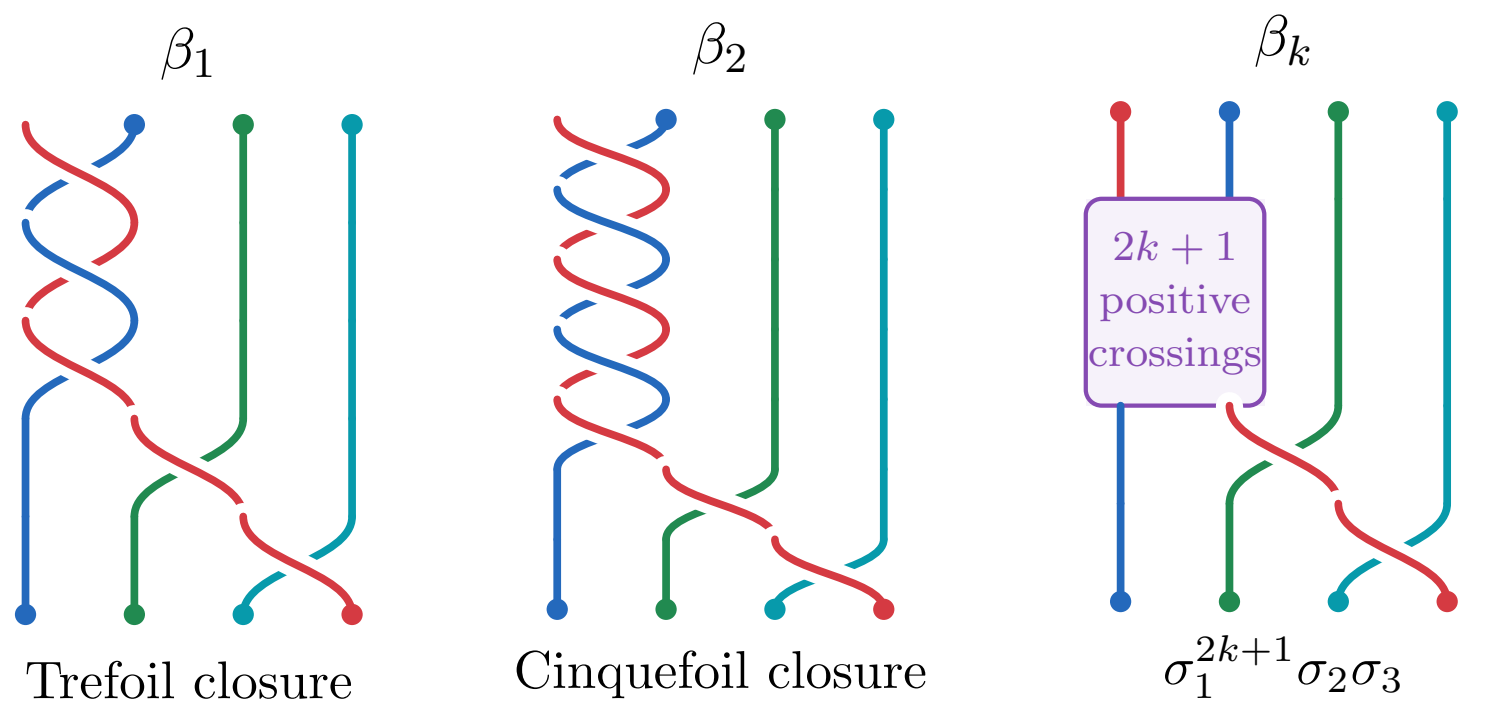
\begin{figure}[htbp]
\centering
\begin{tikzpicture}[x=0.72cm,y=0.72cm,line cap=round,line join=round]
\node at (1.5,0.65) {$\beta_1$};
\fill[strandred] (0,0) circle (2pt);
\fill[strandblue] (1,0) circle (2pt);
\fill[strandgreen] (2,0) circle (2pt);
\fill[strandturquoise] (3,0) circle (2pt);
\draw[strandgreen,line width=1.45pt] (2,0.0000) -- (2,-0.9000);
\draw[strandturquoise,line width=1.45pt] (3,0.0000) -- (3,-0.9000);
\draw[strandblue,line width=1.45pt] (1,0.0000) .. controls (1,-0.4050) and (0,-0.4950) .. (0,-0.9000);
\draw[white,line width=5.1pt] (0,0.0000) .. controls (0,-0.4050) and (1,-0.4950) .. (1,-0.9000);
\draw[strandred,line width=1.45pt] (0,0.0000) .. controls (0,-0.4050) and (1,-0.4950) .. (1,-0.9000);
\draw[strandgreen,line width=1.45pt] (2,-0.9000) -- (2,-1.8000);
\draw[strandturquoise,line width=1.45pt] (3,-0.9000) -- (3,-1.8000);
\draw[strandred,line width=1.45pt] (1,-0.9000) .. controls (1,-1.3050) and (0,-1.3950) .. (0,-1.8000);
\draw[white,line width=5.1pt] (0,-0.9000) .. controls (0,-1.3050) and (1,-1.3950) .. (1,-1.8000);
\draw[strandblue,line width=1.45pt] (0,-0.9000) .. controls (0,-1.3050) and (1,-1.3950) .. (1,-1.8000);
\draw[strandgreen,line width=1.45pt] (2,-1.8000) -- (2,-2.7000);
\draw[strandturquoise,line width=1.45pt] (3,-1.8000) -- (3,-2.7000);
\draw[strandblue,line width=1.45pt] (1,-1.8000) .. controls (1,-2.2050) and (0,-2.2950) .. (0,-2.7000);
\draw[white,line width=5.1pt] (0,-1.8000) .. controls (0,-2.2050) and (1,-2.2950) .. (1,-2.7000);
\draw[strandred,line width=1.45pt] (0,-1.8000) .. controls (0,-2.2050) and (1,-2.2950) .. (1,-2.7000);
\draw[strandblue,line width=1.45pt] (0,-2.7000) -- (0,-3.6000);
\draw[strandturquoise,line width=1.45pt] (3,-2.7000) -- (3,-3.6000);
\draw[strandgreen,line width=1.45pt] (2,-2.7000) .. controls (2,-3.1050) and (1,-3.1950) .. (1,-3.6000);
\draw[white,line width=5.1pt] (1,-2.7000) .. controls (1,-3.1050) and (2,-3.1950) .. (2,-3.6000);
\draw[strandred,line width=1.45pt] (1,-2.7000) .. controls (1,-3.1050) and (2,-3.1950) .. (2,-3.6000);
\draw[strandblue,line width=1.45pt] (0,-3.6000) -- (0,-4.5000);
\draw[strandgreen,line width=1.45pt] (1,-3.6000) -- (1,-4.5000);
\draw[strandturquoise,line width=1.45pt] (3,-3.6000) .. controls (3,-4.0050) and (2,-4.0950) .. (2,-4.5000);
\draw[white,line width=5.1pt] (2,-3.6000) .. controls (2,-4.0050) and (3,-4.0950) .. (3,-4.5000);
\draw[strandred,line width=1.45pt] (2,-3.6000) .. controls (2,-4.0050) and (3,-4.0950) .. (3,-4.5000);
\fill[strandblue] (0,-4.5) circle (2pt);
\fill[strandgreen] (1,-4.5) circle (2pt);
\fill[strandturquoise] (2,-4.5) circle (2pt);
\fill[strandred] (3,-4.5) circle (2pt);
\node[align=center,font=\small] at (1.5,-5.1) {Trefoil closure};
\end{tikzpicture}
\hspace{0.55cm}
\begin{tikzpicture}[x=0.72cm,y=0.72cm,line cap=round,line join=round]
\node at (1.5,0.65) {$\beta_2$};
\fill[strandred] (0,0) circle (2pt);
\fill[strandblue] (1,0) circle (2pt);
\fill[strandgreen] (2,0) circle (2pt);
\fill[strandturquoise] (3,0) circle (2pt);
\draw[strandgreen,line width=1.45pt] (2,0.0000) -- (2,-0.6429);
\draw[strandturquoise,line width=1.45pt] (3,0.0000) -- (3,-0.6429);
\draw[strandblue,line width=1.45pt] (1,0.0000) .. controls (1,-0.2893) and (0,-0.3536) .. (0,-0.6429);
\draw[white,line width=5.1pt] (0,0.0000) .. controls (0,-0.2893) and (1,-0.3536) .. (1,-0.6429);
\draw[strandred,line width=1.45pt] (0,0.0000) .. controls (0,-0.2893) and (1,-0.3536) .. (1,-0.6429);
\draw[strandgreen,line width=1.45pt] (2,-0.6429) -- (2,-1.2857);
\draw[strandturquoise,line width=1.45pt] (3,-0.6429) -- (3,-1.2857);
\draw[strandred,line width=1.45pt] (1,-0.6429) .. controls (1,-0.9321) and (0,-0.9964) .. (0,-1.2857);
\draw[white,line width=5.1pt] (0,-0.6429) .. controls (0,-0.9321) and (1,-0.9964) .. (1,-1.2857);
\draw[strandblue,line width=1.45pt] (0,-0.6429) .. controls (0,-0.9321) and (1,-0.9964) .. (1,-1.2857);
\draw[strandgreen,line width=1.45pt] (2,-1.2857) -- (2,-1.9286);
\draw[strandturquoise,line width=1.45pt] (3,-1.2857) -- (3,-1.9286);
\draw[strandblue,line width=1.45pt] (1,-1.2857) .. controls (1,-1.5750) and (0,-1.6393) .. (0,-1.9286);
\draw[white,line width=5.1pt] (0,-1.2857) .. controls (0,-1.5750) and (1,-1.6393) .. (1,-1.9286);
\draw[strandred,line width=1.45pt] (0,-1.2857) .. controls (0,-1.5750) and (1,-1.6393) .. (1,-1.9286);
\draw[strandgreen,line width=1.45pt] (2,-1.9286) -- (2,-2.5714);
\draw[strandturquoise,line width=1.45pt] (3,-1.9286) -- (3,-2.5714);
\draw[strandred,line width=1.45pt] (1,-1.9286) .. controls (1,-2.2179) and (0,-2.2821) .. (0,-2.5714);
\draw[white,line width=5.1pt] (0,-1.9286) .. controls (0,-2.2179) and (1,-2.2821) .. (1,-2.5714);
\draw[strandblue,line width=1.45pt] (0,-1.9286) .. controls (0,-2.2179) and (1,-2.2821) .. (1,-2.5714);
\draw[strandgreen,line width=1.45pt] (2,-2.5714) -- (2,-3.2143);
\draw[strandturquoise,line width=1.45pt] (3,-2.5714) -- (3,-3.2143);
\draw[strandblue,line width=1.45pt] (1,-2.5714) .. controls (1,-2.8607) and (0,-2.9250) .. (0,-3.2143);
\draw[white,line width=5.1pt] (0,-2.5714) .. controls (0,-2.8607) and (1,-2.9250) .. (1,-3.2143);
\draw[strandred,line width=1.45pt] (0,-2.5714) .. controls (0,-2.8607) and (1,-2.9250) .. (1,-3.2143);
\draw[strandblue,line width=1.45pt] (0,-3.2143) -- (0,-3.8571);
\draw[strandturquoise,line width=1.45pt] (3,-3.2143) -- (3,-3.8571);
\draw[strandgreen,line width=1.45pt] (2,-3.2143) .. controls (2,-3.5036) and (1,-3.5679) .. (1,-3.8571);
\draw[white,line width=5.1pt] (1,-3.2143) .. controls (1,-3.5036) and (2,-3.5679) .. (2,-3.8571);
\draw[strandred,line width=1.45pt] (1,-3.2143) .. controls (1,-3.5036) and (2,-3.5679) .. (2,-3.8571);
\draw[strandblue,line width=1.45pt] (0,-3.8571) -- (0,-4.5000);
\draw[strandgreen,line width=1.45pt] (1,-3.8571) -- (1,-4.5000);
\draw[strandturquoise,line width=1.45pt] (3,-3.8571) .. controls (3,-4.1464) and (2,-4.2107) .. (2,-4.5000);
\draw[white,line width=5.1pt] (2,-3.8571) .. controls (2,-4.1464) and (3,-4.2107) .. (3,-4.5000);
\draw[strandred,line width=1.45pt] (2,-3.8571) .. controls (2,-4.1464) and (3,-4.2107) .. (3,-4.5000);
\fill[strandblue] (0,-4.5) circle (2pt);
\fill[strandgreen] (1,-4.5) circle (2pt);
\fill[strandturquoise] (2,-4.5) circle (2pt);
\fill[strandred] (3,-4.5) circle (2pt);
\node[align=center,font=\small] at (1.5,-5.1) {Cinquefoil closure};
\end{tikzpicture}
\hspace{0.55cm}
\begin{tikzpicture}[x=0.72cm,y=0.72cm,line cap=round,line join=round]
\node at (1.5,0.65) {$\beta_k$};
\fill[strandred] (0,0) circle (2pt);
\draw[strandred,line width=1.45pt] (0,0) -- (0,-0.8);
\fill[strandblue] (1,0) circle (2pt);
\draw[strandblue,line width=1.45pt] (1,0) -- (1,-0.8);
\fill[strandgreen] (2,0) circle (2pt);
\draw[strandgreen,line width=1.45pt] (2,0) -- (2,-0.8);
\fill[strandturquoise] (3,0) circle (2pt);
\draw[strandturquoise,line width=1.45pt] (3,0) -- (3,-0.8);
\fill[twistpurple!7] (-.32,-.8) rectangle (1.32,-2.7);
\draw[twistpurple,line width=.8pt,rounded corners=3pt] (-.32,-.8) rectangle (1.32,-2.7);
\node[align=center,font=\scriptsize,text=twistpurple] at (.5,-1.75) {$2k+1$\\positive\\crossings};
\draw[strandgreen,line width=1.45pt] (2,-.8) -- (2,-2.7);
\draw[strandturquoise,line width=1.45pt] (3,-.8) -- (3,-2.7);
\draw[strandblue,line width=1.45pt] (0,-2.7000) -- (0,-3.6000);
\draw[strandturquoise,line width=1.45pt] (3,-2.7000) -- (3,-3.6000);
\draw[strandgreen,line width=1.45pt] (2,-2.7000) .. controls (2,-3.1000) and (1,-3.2000) .. (1,-3.6000);
\draw[white,line width=5.1pt] (1,-2.7000) .. controls (1,-3.1000) and (2,-3.2000) .. (2,-3.6000);
\draw[strandred,line width=1.45pt] (1,-2.7000) .. controls (1,-3.1000) and (2,-3.2000) .. (2,-3.6000);
\draw[strandblue,line width=1.45pt] (0,-3.6000) -- (0,-4.5000);
\draw[strandgreen,line width=1.45pt] (1,-3.6000) -- (1,-4.5000);
\draw[strandturquoise,line width=1.45pt] (3,-3.6000) .. controls (3,-4.0000) and (2,-4.1000) .. (2,-4.5000);
\draw[white,line width=5.1pt] (2,-3.6000) .. controls (2,-4.0000) and (3,-4.1000) .. (3,-4.5000);
\draw[strandred,line width=1.45pt] (2,-3.6000) .. controls (2,-4.0000) and (3,-4.1000) .. (3,-4.5000);
\fill[strandblue] (0,-4.5) circle (2pt);
\fill[strandgreen] (1,-4.5) circle (2pt);
\fill[strandturquoise] (2,-4.5) circle (2pt);
\fill[strandred] (3,-4.5) circle (2pt);
\node[align=center,font=\small] at (1.5,-5.1) {$\sigma_1^{2k+1}\sigma_2\sigma_3$};
\end{tikzpicture}
\caption{The four-strand braids $\beta_1$, $\beta_2$, and $\beta_k$. The word is read from top to bottom. At each positive crossing the strand entering from the left passes over the strand entering from the right. Colors follow strands through the crossings. The twist box in the third panel represents $2k+1$ positive crossings of the first two strands.}
\label{fig:four-strand-braids}
\end{figure}

The first and third braids belong to the family in Example~\ref{ex:beta_k}
\[
\beta_k=\sigma_1^{2k+1}\sigma_2\sigma_3,\qquad k\geq1,
\]
with $\beta_1=\beta_{3_1}$ and $\beta_2=\beta_{5_1}$.
Its closure is the torus knot $T(2,2k+1)$, whose knot genus is $k$. 
Figure~\ref{fig:four-strand-braids} shows $\beta_1$, $\beta_2$, and
$\beta_k$.

Applying
Proposition~\ref{prop:braid-surface} over $\Sigma_2$ gives a connected
symplectic multisection $F_{\beta_k}\subset\Sigma_2\times S^2$ with
\[
[F_{\beta_k}]=4[\Sigma_2\times\{\mathrm{pt}\}],\qquad
[F_{\beta_k}]^2=0,\qquad
g(F_{\beta_k})=4(2-1)+1=5.
\]
Thus the genus of the multisection is independent of the knot genus.
Moreover, the coverings $F_{\beta_k}\to\Sigma_2$ are equivalent for
all $k$: their monodromy sends $a_1$ to the same four-cycle and all
other standard generators to the identity. The braid action on the
punctured sphere can nevertheless retain information beyond this
permutation, as explained in Section~\ref{sec:braids}. The peripheral
fundamental-group maps are used in Section~\ref{sec:fixed-gluing} to
distinguish sums for a fixed gluing.

The companion surface is the connected cyclic degree-two multisection
$F_{3,2}\subset\Sigma_3\times S^2$ of
Proposition~\ref{prop:braid-surface}. It satisfies
\[
[F_{3,2}]=2[\Sigma_3\times\{\mathrm{pt}\}],\qquad
[F_{3,2}]^2=0,\qquad g(F_{3,2})=2(3-1)+1=5.
\]
Thus $F_{\beta_k}$ and $F_{3,2}$ have the same genus and trivial normal
bundles. The following proposition describes the resulting symplectic
sums and a gluing for which the first Betti number vanishes.

\begin{proposition}\label{prop:four-sum}
Let $F_{3,2}\subset\Sigma_3\times S^2$ be the connected cyclic
degree-two multisection of Proposition~{\rm \ref{prop:braid-surface}}. 
For every $k\geq 1$, the surfaces $F_{\beta_k}$ and $F_{3,2}$ admit a symplectic sum. Every
such sum is aspherical and minimal, and has Euler characteristic $4$
and signature $0$. The surface gluing can be chosen so that the first
Betti number is zero.
\end{proposition}

\begin{proof}
Both surfaces have genus five and square zero. Rescale one of the
summand symplectic forms so that the surfaces have equal area. For
any orientation-preserving diffeomorphism between them, Moser's
theorem gives an isotopic diffeomorphism identifying their area forms.
The trivial normal bundles can then be identified with opposite
normal-circle orientations, and Gompf's construction gives the
symplectic sum as in~\cite{Gompf}. The sum formulas give
\[
 \begin{aligned}
 e(Z)
 &=e(\Sigma_2\times S^2)+e(\Sigma_3\times S^2)-2e(\Sigma_5)\\
 &=-4-8-2(-8)=4,\\
 \sigma(Z)
 &=\sigma(\Sigma_2\times S^2)+\sigma(\Sigma_3\times S^2)=0.
 \end{aligned}
\]
Asphericity and minimality follow from Proposition~\ref{prop:aspherical-sums}.

We give an explicit surface gluing with $b_1(Z)=0$. Let
$q_{\beta_k}:F_{\beta_k}\to\Sigma_2$ and
$q_{3,2}:F_{3,2}\to\Sigma_3$ be the covering maps, and write
\[
 \begin{aligned}
 Q_{4,2}=(q_{\beta_k})_* &: 
 H_1(F_{\beta_k};\Z)\longrightarrow H_1(\Sigma_2;\Z),\\
 Q_{2,3}=(q_{3,2})_* &: 
 H_1(F_{3,2};\Z)\longrightarrow H_1(\Sigma_3;\Z).
 \end{aligned}
\]
Choose symplectic bases $(a,b,c,d)$ and
$(\alpha,\eta,\gamma,\delta,\epsilon,\zeta)$ on the two base
surfaces, with covering monodromy supported on $a$ and $\alpha$,
respectively. Cutting the bases along $b$ and $\eta$ and gluing
four or two copies cyclically gives symplectic bases
$(e_j,f_j)_{j=0}^4$ and $(e'_j,f'_j)_{j=0}^4$ upstairs with
\[
 \begin{array}{c|ccccc}
 j&0&1&2&3&4\\ \hline
 Q_{4,2}(e_j)&4a&c&c&c&c\\
 Q_{4,2}(f_j)&b&d&d&d&d\\
 Q_{2,3}(e'_j)&2\alpha&\gamma&\epsilon&\gamma&\epsilon\\
 Q_{2,3}(f'_j)&\eta&\delta&\zeta&\delta&\zeta
 \end{array}
\]
Here the first pairs $(e_0,f_0)$ and $(e'_0,f'_0)$ come from the cyclic cover of the first handles containing 
$(a,b)$ and $(\alpha,\eta)$, respectively;
the other pairs are the lifts of the remaining handles, ordered
sheet by sheet.

Define an integral symplectic isomorphism $S$ by
\[
 \begin{aligned}
 S(e_0,f_0)&=(e'_1,f'_1),&
 S(e_1,f_1)&=(e'_0,f'_0),\\
 S(e_j,f_j)&=(e'_j,f'_j)\quad(j=2,3),&
 S(e_4,f_4)&=(-e'_4,-f'_4).
 \end{aligned}
\]
This exchanges the first two symplectic pairs and changes the signs
of both elements of the last pair, so it preserves the intersection
form. It is induced by an orientation-preserving diffeomorphism
$\psi:F_{\beta_k}\to F_{3,2}$, by the surjectivity of the symplectic
representation (cf.~\cite[Theorem~6.4]{FarbMargalit}). As above, $\psi$ can be chosen to
identify the area forms without changing its action on homology.

By Lemma~\ref{lem:h1}, projection identifies the rational first homology of
each exterior with that of its base. The meridians vanish over
$\Q$, and the boundary surface maps are therefore $Q_{4,2}$
and $Q_{2,3}$. Since the boundary and the two exteriors are connected,
Mayer--Vietoris gives
\[
 H_1(Z;\Q)\cong
 \frac{H_1(\Sigma_2;\Q)\oplus H_1(\Sigma_3;\Q)}
 {\{(Q_{4,2}v,-Q_{2,3}Sv)\mid v\in H_1(F_{\beta_k};\Q)\}}.
\]
The relations in this quotient are
\[
 \begin{aligned}
 4a&=\gamma,& b&=\delta,\\
 c&=2\alpha,& d&=\eta,\\
 c&=\epsilon,& d&=\zeta,\\
 c&=\gamma,& d&=\delta,\\
 c&=-\epsilon,& d&=-\zeta.
 \end{aligned}
\]
The third and fifth rows give $2c=2d=0$, hence $c=d=0$ over
$\Q$. The remaining relations then kill all the other
generators. Thus $H_1(Z;\Q)=0$. The covering maps are
independent of $k$, so the same homological
gluing works for every $k$.
\end{proof}

Finally, the double cover of $S^2$ branched at the four braid points
gives a useful invariant of the spherical braid monodromy. This
double cover is a torus. A half-twist about an arc joining two branch
points lifts to a Dehn twist about its inverse image. In a suitable
symplectic basis, the lifts of the three adjacent half-twists act by
\[
\sigma_1,\sigma_3\longmapsto
A=\begin{pmatrix}1&1\\0&1\end{pmatrix},
\qquad
\sigma_2\longmapsto
B=\begin{pmatrix}1&0\\-1&1\end{pmatrix}.
\]
Using the composition convention of Section~\ref{sec:braids}, we obtain
\[
 A^{2k+1}BA
 =\begin{pmatrix}1&2k+1\\0&1\end{pmatrix}
  \begin{pmatrix}1&0\\-1&1\end{pmatrix}
  \begin{pmatrix}1&1\\0&1\end{pmatrix}
 =\begin{pmatrix}-2k&1\\-1&0\end{pmatrix}.
\]
The two lifts of a spherical diffeomorphism differ by the elliptic
involution, which acts as $-I$ on the first homology of the torus. 
Changing the symplectic basis conjugates the matrix.  
Consequently, the absolute value of the trace is a well-defined conjugacy invariant of the spherical mapping
class. It equals $2k$ for $\beta_k$, so these spherical mapping
classes are pairwise nonconjugate. 
They remain distinct even if inversion of the monodromy is
allowed, since a matrix in $\mathrm{SL}(2,\mathbb Z)$ and
its inverse have the same trace.

The torus in this calculation is the double cover of the sphere
fiber branched at the four points of the stabilized braid. For
example, its monodromy for $\beta_1$ has absolute trace $2$, whereas
the capped trefoil-fiber monodromy at the beginning of this section
has absolute trace $1$. As explained in Section~\ref{sec:braids}, this invariant
distinguishes the marked sphere monodromies; it does not establish
ambient non-isotopy of the multisections or distinguish the closed
sums for a fixed gluing. The latter distinction is obtained in
Section~\ref{sec:fixed-gluing} from the peripheral fundamental-group maps.

\section{Five fake-quadric patterns}
\label{sec:five-patterns}

Write $X_g=\Sigma_g\times S^2$, and let $E_g\to\Sigma_g$ be the
nontrivial oriented sphere bundle.  We use the five pairs considered in
\cite{AY}:
\[
\begin{array}{c|c|c|c}
\text{case} & \text{first surface} & \text{second surface} & G\\ \hline
\text{product }(2,3) & F_{4,2}\subset X_4 & F_{3,3}\subset X_3 & 7\\
\text{product }(2,4) & F_{3,2}\subset X_3 & F_{2,4}\subset X_2 & 5\\
\text{product }(3,3) & F_{2,3}\subset X_2 & F_{2,3}\subset X_2 & 4\\
\text{twisted }(2,4) & A_{3,2}\subset E_3 & A_{2,4}\subset E_2 & 5\\
\text{hybrid }(2,4) & F_{3,2}\subset X_3 & A_{2,4}\subset E_2 & 5
\end{array}
\]
The pair $(p,q)$ records the degrees of the first and second
multisections, and $G$ is their common genus.  Here $F_{g,p}$ is the
cyclic product multisection with
\[
 [F_{g,p}]=p[\Sigma_g\times\{\mathrm{pt}\}],
\]
while the twisted multisections have
\[
 [A_{g,p}]=p[S_g]-\frac p2[D],\qquad
 [S_g]^2=1,\quad [D]^2=0,\quad [S_g]\cdot[D]=1,
\]
where $S_g$ is a section, $D$ is a sphere fiber, and $p$ is even.
Both multisection types have self-intersection zero.  

Writing $g,h$ for the two base genera,
the unbranched covering formula gives the common genus
\[
 G=p(g-1)+1=q(h-1)+1.
\]
To obtain a sum $Z$\/ with Euler characteristic four, we also require
$G=g+h$.  Indeed, the symplectic sum formula gives
\[
 \begin{aligned}
 e(Z)&=4(1-g)+4(1-h)-2(2-2G)\\
     &=4(1+G-g-h),
 \end{aligned}
\]
so $e(Z)=4$ if and only if $G=g+h$.
The five pairs in the above table satisfy this condition.  
Thus every sum along one of these pairs has
$e(Z)=4$.  It also has $\sigma(Z)=0$, by additivity and the vanishing
of the signatures of the two sphere bundles.  The calculations below
give product gluings with finite first homology; the twisted and hybrid
gluings are treated in Section~\ref{sec:five-braids}.

We next specify the homological effect of the normal framing.
Unless otherwise indicated, all homology groups below have integer
coefficients. Let $M=E\setminus\Int\nu F$ be the exterior of a
connected square-zero multisection $q:F\to\Sigma_g$ in an oriented
sphere bundle $E\to\Sigma_g$. The normal bundle of $F$ is an
oriented real rank-two bundle. Excision and the homological Thom
isomorphism therefore give
\[
 H_k(E,M)\cong H_k(\nu F,\partial\nu F)\cong H_{k-2}(F).
\]
In the present square-zero setting, the normal bundle is trivial, so
this identification also follows directly from the relative K\"unneth
formula for $(F\times D^2,F\times S^1)$. Since $F$ is connected,
\[
 H_2(E,M)\cong H_0(F)\cong\Z,\qquad H_1(E,M)=0.
\]
The positive generator of $H_2(E,M)$ is represented by an oriented
normal disk $\Delta$ meeting $F$ once positively, with
$\partial\Delta\subset\partial\nu F\subset M$.

The long exact sequence of the pair $(E,M)$ contains
\[
 H_2(E)\xrightarrow{j_*}H_2(E,M)
 \xrightarrow{\partial}H_1(M)
 \longrightarrow H_1(E)\longrightarrow H_1(E,M).
\]
The connecting homomorphism sends a relative cycle to the homology
class of its boundary. Thus $\partial[\Delta]=[\partial\Delta]=\mu$,
where $\mu$ is the positively oriented meridian class.
Under the Thom identification, $j_*$ is intersection with $[F]$:
a surface transverse to $F$ contributes one normal disk for each
intersection point, counted with sign. Consequently, the sequence
becomes
\[
 H_2(E)\xrightarrow{\,\cdot[F]\,}\Z
 \xrightarrow{\,1\mapsto\mu\,}H_1(M)
 \longrightarrow H_1(E)\longrightarrow0.
\]

If the intersection map has image $d\Z$, exactness says that the
kernel of $n\mapsto n\mu$ is precisely $d\Z$. Hence $\mu$ has
order exactly $d$. The sphere-bundle projection identifies
$H_1(E)$ with $H_1(\Sigma_g)$, so we obtain
\[
 0\longrightarrow\Z/d
 \xrightarrow{\,\overline{1}\mapsto\mu\,}H_1(M)
 \longrightarrow H_1(\Sigma_g)\longrightarrow0.
\]
Choosing lifts of a basis of the free abelian group $H_1(\Sigma_g)$
splits this sequence, giving a noncanonical isomorphism
\[
 H_1(M)\cong H_1(\Sigma_g)\oplus\Z/d.
\]
For the product class $[F_{g,p}]$, intersection with the section and
fiber gives $0$ and $p$, so $d=p$. For the twisted class
$[A_{g,p}]=p[S_g]-(p/2)[D]$, these intersection numbers are $p/2$
and $p$, so $d=p/2$.

\begin{lemma}\label{lem:peripheral}
Let $F_1,F_2\subset E$ be connected square-zero multisections of an
oriented sphere bundle $E\to\Sigma_g$, with covering maps
$q_i:F_i\to\Sigma_g$. Put
\[
 M_i=E\setminus\Int\nu F_i,\qquad i=1,2,
\]
and let $\mu_i\in H_1(M_i)$ be the positive meridian classes.
Suppose that both $\mu_i$ have order $d$ and that there is a
diffeomorphism $h:F_1\to F_2$ satisfying $q_2\circ h=q_1$.
Fix isomorphisms
\[
 \Theta_i:H_1(M_i)\longrightarrow H_1(\Sigma_g)\oplus\Z/d
\]
whose first components are induced by the bundle projection and
which satisfy $\Theta_i(\mu_i)=(0,\overline{1})$.
Then there are oriented normal framings for which the boundary
parametrizations
\[
 \iota_i:F_i\times S^1\longrightarrow\partial M_i\subset M_i
\]
satisfy
\[
 \Theta_i\bigl((\iota_i)_*(\xi,n)\bigr)
   =\bigl((q_i)_*\xi,\overline n\bigr),
 \qquad \xi\in H_1(F_i),\quad n\in\Z.
\]
In particular,
\[
 \Theta_2\circ(\iota_2)_*\circ
 (h_*\oplus\mathrm{id}_{\Z})
 =\Theta_1\circ(\iota_1)_*.
\]
Consequently, if the framed gluing maps correspond under
$h\times\mathrm{id}_{S^1}$, the above identifications carry the
integral Mayer--Vietoris relations for the sum along $F_1$ to
those for the sum along $F_2$. In particular, the two sums have
isomorphic first homology groups.
\end{lemma}

\begin{proof}
Fix $i\in\{1,2\}$. Since $F_i^2=0$, its oriented normal bundle is
trivial. Choose an initial normal framing and write
$\iota_i^0:F_i\times S^1\to M_i$ for the corresponding boundary
parametrization followed by inclusion. A pushed-off loop in $F_i$
projects to a loop homologous to its image under $q_i$, while the
positive normal circle represents $\mu_i$. Hence
\[
 \Theta_i\bigl((\iota_i^0)_*(\xi,n)\bigr)
 =\bigl((q_i)_*\xi,\lambda_i(\xi)+\overline n\bigr)
\]
for a homomorphism $\lambda_i:H_1(F_i)\to\Z/d$.

Because $H_1(F_i)$ is free abelian, there is a homomorphism
$\widetilde\lambda_i:H_1(F_i)\to\Z$ whose reduction modulo $d$
is $-\lambda_i$. Choose a smooth map $f_i:F_i\to S^1$ inducing
$\widetilde\lambda_i$ on the first homology, and change the framing by
setting
\[
 \iota_i(x,z)=\iota_i^0(x,f_i(x)z).
\]
The map $(x,z)\mapsto(x,f_i(x)z)$ induces
$(\xi,n)\mapsto(\xi,n+\widetilde\lambda_i(\xi))$ on the first homology.
Therefore
\[
 \begin{aligned}
 \Theta_i\bigl((\iota_i)_*(\xi,n)\bigr)
 &=\bigl((q_i)_*\xi,
   \lambda_i(\xi)+\overline n+
   \overline{\widetilde\lambda_i(\xi)}\bigr)\\
 &=\bigl((q_i)_*\xi,\overline n\bigr).
 \end{aligned}
\]
We call this framing \emph{adapted} to the chosen splitting
$\Theta_i$.

Apply this construction for $i=1,2$. Since
$(q_2)_*\circ h_*=(q_1)_*$, the two boundary maps satisfy the asserted
identity. Thus, under $h_*\oplus\mathrm{id}_{\Z}$ on the boundary
and $\Theta_2^{-1}\circ\Theta_1$ on the exterior homology, their
inclusion maps agree. Transporting the gluing to any fixed second
exterior through $h\times\mathrm{id}_{S^1}$ therefore identifies the
two Mayer--Vietoris relation maps, proving the final assertion.
\end{proof}

An adapted framing need not be unique: changing it by an element of
$dH^1(F;\Z)$ leaves the boundary map on first homology
unchanged.  In particular, the lemma identifies only the homological
boundary maps.  It does not identify the peripheral fundamental-group
maps or the boundary spin structures.

\subsection{The homology matrices}
\label{subsec:matrices}

For the cyclic degree-$p$ covering determined by
$a_1\mapsto 1\in\Z/p$ and all other standard base generators
mapping to zero, put $G=p(g-1)+1$.  Cutting the base along a curve
representing $b_1$ and cyclically joining $p$ copies of the resulting
surface gives a symplectic basis upstairs consisting of
\[
 (e_0,f_0),\qquad (e_{i,r},f_{i,r}),
 \quad 1\leq i\leq g-1,\quad 0\leq r<p.
\]
Here $e_0$ traverses the $p$ sheets over $a_1$, $f_0$ is a lift of
$b_1$, and the remaining pairs are the lifts of the other base
handles.  In the base basis $(a_1,b_1,\ldots,a_g,b_g)$, the covering
pushforward
\[
 Q_{p,g}:H_1(\Sigma_G)\longrightarrow H_1(\Sigma_g)
\]
is therefore given by
\[
 \begin{aligned}
 Q_{p,g}(e_0)&=pa_1,&\quad Q_{p,g}(f_0)&=b_1,\\
 Q_{p,g}(e_{i,r})&=a_{i+1},&
 Q_{p,g}(f_{i,r})&=b_{i+1}.
 \end{aligned}
\]
Replacing the cyclic braid on $a_1$ by a braid with the same strand
permutation leaves the covering unchanged and hence leaves
$Q_{p,g}$ unchanged in these bases.  The peripheral
fundamental-group maps may nevertheless change.

Choose adapted framings on both exteriors and glue their boundaries by
\[
 (x,z)\longmapsto(\psi(x),\overline z),
\]
so that the normal twist is zero.  Let
$S=\psi_*\in\operatorname{Sp}(2G,\Z)$, with $\psi$ mapping the
\emph{first} surface in the table to the \emph{second}.  More generally,
write $Q_1,Q_2$ for the two covering pushforwards and $d_1,d_2$ for
the meridian orders.  

Let $M_1,M_2$ be the two exteriors and put
$B=\partial M_1\cong\Sigma_G\times S^1$.
Since $B,M_1,M_2$ are connected, the Mayer--Vietoris sequence gives
\[
H_1(B)\xrightarrow{\alpha}H_1(M_1)\oplus H_1(M_2)
\longrightarrow H_1(Z)\longrightarrow\Z
\xrightarrow{\,n\mapsto(n,-n)\,}\Z^2.
\]
The last map is injective, so $H_1(Z)\cong\operatorname{coker}\alpha$.
In the adapted framings, with zero normal twist, the gluing reverses
the normal circle. Thus, grouping the base and meridian summands,
\[
\alpha(\xi,n)
=\bigl((Q_1\xi,-Q_2S\xi),(\overline n,\overline n)\bigr).
\]
The surface and meridian relations therefore separate.
The former give the cokernel of
$\binom{Q_1}{-Q_2S}$, while the latter give
\[
\frac{\Z/d_1\oplus\Z/d_2}
{\langle(\overline1,\overline1)\rangle}
\cong\Z/\gcd(d_1,d_2).
\]
Indeed, the meridian relations are
$d_1t=0$, $d_2s=0$, and $t+s=0$; eliminating $s$ gives the
displayed cyclic group. Consequently,
\begin{equation}\label{eq:homology-formula}
H_1(Z)\cong
\operatorname{coker}
\begin{pmatrix}Q_1\\-Q_2S\end{pmatrix}
\oplus\Z/\gcd(d_1,d_2).
\end{equation}
For a product pair, $Q_1=Q_{p,g}$, $Q_2=Q_{q,h}$, and
$(d_1,d_2)=(p,q)$.

We now give explicit gluings for the three product rows.  For
$c\in\Z$ and $m\geq1$, set
\[
 U_c=\begin{pmatrix}1&c\\0&1\end{pmatrix},\qquad
 L_c=\begin{pmatrix}1&0\\c&1\end{pmatrix},\qquad
 R_m=\begin{pmatrix}2-m&-1\\m-1&1\end{pmatrix},\qquad
 T_m=\begin{pmatrix}1&m-1\\1&m\end{pmatrix}.
\]
Each block has determinant one.  Number the symplectic pairs from
zero, starting with $(e_0,f_0)$.  In the $(2,3)$ and $(3,3)$ rows,
order the remaining pairs lexicographically by $(i,r)$ on each side.
For example, in the $(2,3)$ case, the source pairs numbered
$1,2$ are $(e_{1,0},f_{1,0})$ and $(e_{1,1},f_{1,1})$.
On the target surface, the pairs numbered $1,2,3$ are
$(e'_{1,0},f'_{1,0})$, $(e'_{1,1},f'_{1,1})$, and
$(e'_{1,2},f'_{1,2})$, respectively.
In the $(2,4)$ row, order them by $(r,i)$, which is 
the ordering of the peripheral generators used in
Section~\ref{sec:fixed-gluing}.
For each $m\geq1$, we specify a matrix
$S_m\in\operatorname{Sp}(2G,\Z)$ representing the induced map
on first homology of a surface gluing
$\psi_m:F_1\to F_2$, in the ordered symplectic bases chosen above.
Here the integer $m\geq1$ parametrizes the choice of gluing.
For each product pair, varying $m$ will give sums with
pairwise nonisomorphic first homology groups.

To describe $S_m$, write the ordered symplectic pairs on the
first surface as $(e_j,f_j)$ and those on the second surface as
$(e'_j,f'_j)$, where $0\leq j<G$. Here we number the pairs
in the order specified above.

The map $S_m$ sends the $j$-th source pair to the $\pi(j)$-th
target pair according to a matrix $B_j$. Explicitly, if
\[
B_j=\begin{pmatrix}a&b\\c&d\end{pmatrix},
\]
then
\[
\begin{aligned}
S_m(e_j)&=a e'_{\pi(j)}+c f'_{\pi(j)},\\
S_m(f_j)&=b e'_{\pi(j)}+d f'_{\pi(j)}.
\end{aligned}
\]
The following table specifies these maps. An entry $(\ell,B)$
in row $j$ means that $\pi(j)=\ell$ and $B_j=B$.
We write $I$ for the $2\times2$ identity matrix; a dash indicates
that the source pair does not occur in that case.
\[
\begin{array}{c|c|c|c}
\text{source pair }j
 &(p,q)=(2,3)&(p,q)=(2,4)&(p,q)=(3,3)\\ \hline
0 &(6,U_1)&(3,L_1)&(3,U_1)\\
1 &(1,L_m)&(1,R_m)&(0,U_1)\\
2 &(3,U_1)&(2,L_{-1})&(2,U_{-1})\\
3 &(0,I)&(4,L_{-1})&(1,T_m)\\
4 &(4,U_1)&(0,U_{-1})&\text{--}\\
5 &(5,I)&\text{--}&\text{--}\\
6 &(2,L_1)&\text{--}&\text{--}
\end{array}
\]
For example, in the $(2,4)$ case, the entry $(3,L_1)$ in row
zero gives
\[
S_m(e_0)=e'_3+f'_3,\qquad S_m(f_0)=f'_3.
\]

In each case, every target pair occurs exactly once, and every
matrix $B_j$ has determinant one. Thus $S_m$ preserves the
intersection pairing within each symplectic pair, while distinct
source pairs map to distinct target pairs and remain orthogonal.
Hence $S_m\in\operatorname{Sp}(2G,\Z)$.

The surjectivity of the symplectic representation of the mapping class group
\cite[Theorem~6.4]{FarbMargalit} realizes each $S_m$ by an orientation-preserving
surface gluing.  After rescaling the symplectic forms to give the
summing surfaces equal area, Moser's theorem gives an area-preserving
representative in the same isotopy class, as in Section~\ref{sec:four-strand}.
In the $(2,4)$ row, $S_m$ maps the degree-two surface to the
degree-four surface.  The gluing in Section~\ref{sec:fixed-gluing} is written in the
opposite direction and consequently has matrix $S_m^{-1}$.

For each product pair, fix an area-preserving gluing
$\psi_m:F_{g,p}\to F_{h,q}$ inducing $S_m$, and let
\[
Z_m=
\bigl(X_g\setminus\Int\nu F_{g,p}\bigr)
\cup_{\varphi_m}
\bigl(X_h\setminus\Int\nu F_{h,q}\bigr),
\qquad
\varphi_m(x,z)=(\psi_m(x),\overline z),
\]
where the boundary identifications use the adapted normal
framings. Thus $Z_m$ is the resulting symplectic sum, with
$m$ recording the choice of gluing. We use this notation
separately for each product pair.

To compute the cokernels explicitly, we distinguish the homology
bases of the summing surfaces from those of the base surfaces.
We retain $(e_j,f_j)$ and $(e'_j,f'_j)$ for the two summing
surfaces, and write
\[
(a_1,b_1,\ldots,a_g,b_g),\qquad
(a'_1,b'_1,\ldots,a'_h,b'_h)
\]
for the bases of $H_1(\Sigma_g)$ and $H_1(\Sigma_h)$, respectively.
The map
\[
\xi\longmapsto
\bigl(Q_{p,g}\xi,-Q_{q,h}S_m\xi\bigr)
\]
from the homology of the first summing surface to
$H_1(\Sigma_g)\oplus H_1(\Sigma_h)$ has matrix
\[
M_m=\begin{pmatrix}Q_{p,g}\\-Q_{q,h}S_m\end{pmatrix}.
\]
Its rows correspond to the base generators listed above, and
its columns correspond to the ordered basis elements of the
first summing surface. Each column therefore gives a relation
among the base generators in $\operatorname{coker}M_m$.
Since $G=g+h$, the matrix $M_m$ is square of size $2G$.

\smallskip
\noindent\emph{The product $(2,3)$ case.}
Here $Q_1=Q_{2,4}$ and $Q_2=Q_{3,3}$. Each source basis
element $\xi$ gives the relation
$Q_{2,4}(\xi)=Q_{3,3}(S_m\xi)$ in $\operatorname{coker}M_m$.
All equalities below are understood in this cokernel.

The source pairs numbered $1,2$ both project to $(a_2,b_2)$.
They map to the target pairs numbered $1,3$, respectively,
both of which project to $(a'_2,b'_2)$. The corresponding
blocks $L_m$ and $U_1$ therefore give
\[
\begin{array}{c|cc}
\text{source pair}&\multicolumn{2}{c}{\text{relations}}\\ \hline
1 & a_2=a'_2+mb'_2 & b_2=b'_2\\
2 & a_2=a'_2       & b_2=a'_2+b'_2
\end{array}
\]
Comparing the two expressions for $b_2$ gives $a'_2=0$.
The relations for $a_2$ then give $a_2=0$ and $mb'_2=0$.

The remaining source pairs give
\[
\begin{array}{c|cc}
\text{source pair}&\multicolumn{2}{c}{\text{relations}}\\ \hline
0 & 2a_1=a'_3       & b_1=a'_3+b'_3\\
3 & a_3=3a'_1       & b_3=b'_1\\
4 & a_3=a'_3        & b_3=a'_3+b'_3\\
5 & a_4=a'_3        & b_4=b'_3\\
6 & a_4=a'_2+b'_2   & b_4=b'_2
\end{array}
\]
Using $a'_2=0$, the relations from pairs $0,5,6$ imply
\[
a'_3=a_4=b'_2=b'_3=b_4=2a_1,
\qquad b_1=4a_1.
\]
The relations from pairs $3,4$ then give
\[
a_3=2a_1=3a'_1,\qquad b_3=b'_1=4a_1.
\]
Together with $b_2=b'_2$, these eliminate every generator
except $a_1,a'_1$. The remaining relations are
$2a_1=3a'_1$ and $2ma_1=0$, the latter coming from $mb'_2=0$.
Thus
\[
\operatorname{coker}M_m
=\langle a_1,a'_1\mid 2a_1-3a'_1=0,\;2ma_1=0\rangle.
\]
The first relation alone gives an infinite cyclic group
generated by $a_1-a'_1$, with
\[
a_1=3(a_1-a'_1),\qquad a'_1=2(a_1-a'_1).
\]
The second relation becomes $6m(a_1-a'_1)=0$, so
$\operatorname{coker}M_m\cong\Z/(6m)$.

\smallskip
\noindent\emph{The product $(2,4)$ case.}
Here $Q_1=Q_{2,3}$ and $Q_2=Q_{4,2}$. With the ordering
by $(r,i)$, the source pairs numbered $1,3$ both project
to $(a_2,b_2)$. Their target pairs, numbered $1,4$,
both project to $(a'_2,b'_2)$. The blocks $R_m$ and
$L_{-1}$ therefore give the following relations in
$\operatorname{coker}M_m$:
\[
\begin{aligned}
a_2&=(2-m)a'_2+(m-1)b'_2,
& b_2&=-a'_2+b'_2,\\
a_2&=a'_2-b'_2,
& b_2&=b'_2.
\end{aligned}
\]
Comparing the expressions for $b_2$ gives $a'_2=0$.
The relations for $a_2$ then give
$a_2=-b'_2$ and $mb'_2=0$.

The source pairs numbered $0,2,4$ give, respectively,
\[
\begin{aligned}
2a_1&=a'_2+b'_2,
& b_1&=b'_2,\\
a_3&=a'_2-b'_2,
& b_3&=b'_2,\\
a_3&=4a'_1,
& b_3&=-4a'_1+b'_1.
\end{aligned}
\]
Since $a'_2=0$, these imply
\[
b'_2=2a_1,\qquad a_3=-2a_1=4a'_1,\qquad b'_1=0.
\]
Every other generator is now determined by $a_1,a'_1$,
and the remaining relations are
\[
\operatorname{coker}M_m
=\langle a_1,a'_1\mid 2a_1+4a'_1=0,\;2ma_1=0\rangle.
\]
Use the new basis $(a_1+2a'_1,a'_1)$ and subtract $m$
times the first relation from the second. The relations become
\[
2(a_1+2a'_1)=0,\qquad 4ma'_1=0.
\]
Hence
\[
\operatorname{coker}M_m\cong\Z/2\oplus\Z/(4m).
\]

\smallskip
\noindent\emph{The product $(3,3)$ case.}
Both coverings have degree three over a genus-two surface.
On each side, the pairs numbered $1,2,3$ are the three
lifts of the second base handle. Thus they project to
$(a_2,b_2)$ on the source side and $(a'_2,b'_2)$ on
the target side.

The source pairs numbered $2,3$ map to the target pairs
numbered $2,1$, with blocks $U_{-1}$ and $T_m$, respectively.
They give the following relations in $\operatorname{coker}M_m$:
\[
\begin{aligned}
a_2&=a'_2,
& b_2&=-a'_2+b'_2,\\
a_2&=a'_2+b'_2,
& b_2&=(m-1)a'_2+mb'_2.
\end{aligned}
\]
Comparing the expressions for $a_2$ gives $b'_2=0$.
The relations for $b_2$ then give
\[
b_2=-a'_2,\qquad ma'_2=0.
\]

The source pairs numbered $0,1$ both have block $U_1$
and map to the target pairs numbered $3,0$, respectively.
Their relations are
\[
\begin{aligned}
3a_1&=a'_2,
& b_1&=a'_2+b'_2,\\
a_2&=3a'_1,
& b_2&=3a'_1+b'_1.
\end{aligned}
\]
Together with the preceding relations, these imply
\[
a_2=a'_2=3a_1=3a'_1,\qquad
b_1=3a_1,\qquad b_2=-3a_1,\qquad b'_1=-6a_1.
\]
Hence every generator is determined by $a_1,a'_1$, and
the remaining relations are $3a_1=3a'_1$ and $3ma_1=0$,
the latter coming from $ma'_2=0$. Therefore
\[
\operatorname{coker}M_m
=\langle a_1,a'_1\mid 3a_1-3a'_1=0,\;3ma_1=0\rangle.
\]
In the basis $(a_1-a'_1,a_1)$, these are the independent
relations
\[
3(a_1-a'_1)=0,\qquad 3ma_1=0.
\]
It follows that
\[
\operatorname{coker}M_m\cong\Z/3\oplus\Z/(3m).
\]

\smallskip
The following table summarizes the first homology groups of $Z_m$, 
including the meridian summands from (\ref{eq:homology-formula}).
\[
\begin{array}{c|c}
 (p,q)&H_1(Z_m;\Z)\\ \hline
 (2,3)&\Z/(6m)\\
 (2,4)&(\Z/2)^2\oplus\Z/(4m)\\
 (3,3)&(\Z/3)^2\oplus\Z/(3m)
\end{array}
\]
In particular, the orders of the first homology groups are $6m$, $16m$, and $27m$.

\section{Genus-two knots with two three-point orbits}
\label{sec:genus-two}

We next obtain degree-three product multisections from a family
of genus-two fibered knots. Let
\[
J_0=3_1\mathbin{\#}\overline{3_1}
\]
be the square knot, and let $J_n$ be the knot obtained by the
$n$ Stallings twists in the construction of Park and Yun in~\cite{ParkYun}.
These knots are fibered and have genus two.  
They all have the same Alexander polynomial $(t^2-t+1)^2$.  
We first distinguish the knots and then determine the branch-point permutation of the hyperelliptic quotient for the even-indexed subfamily.

By \cite[\S3, p.~565]{QuachWeber} and \cite[preprint version, Lemma~3.2(3), p.~7]{ParkYun},
the Alexander module of $J_n$ has a presentation over
$\Lambda=\mathbb Z[t,t^{-1}]$ with matrix
\[
\begin{pmatrix}
f(t)&nt\\
0&f(t)
\end{pmatrix},
\qquad
f(t)=t^2-t+1.
\]
Let $\mathcal A_n$ denote the Alexander module of $J_n$.
For $n>0$, consider the quotient
$\mathcal A_n/f(t)\mathcal A_n$. Reducing the displayed
presentation modulo $f(t)$ gives the matrix
\[
\begin{pmatrix}
0&nt\\
0&0
\end{pmatrix}
\]
over $\Lambda/(f(t))$. Since $t$ is a unit, we obtain
\[
\mathcal A_n/f(t)\mathcal A_n
\cong
\Lambda/(f(t))\oplus\Lambda/(f(t),n).
\]
The identity $t(1-t)=1$ modulo $f(t)$ allows us to eliminate
negative powers of $t$ (by $
t^{-k}=(1-t)^k$).  
Moreover, since $f(t)$ is monic of
degree two, every element of these two summands has a unique
representative $a+bt$, with coefficients in $\mathbb Z$ and
$\mathbb Z/n$, respectively. Consequently, as abelian groups,
\[
\mathcal A_n/f(t)\mathcal A_n
\cong
\mathbb Z^2\oplus(\mathbb Z/n)^2.
\]
Its torsion subgroup therefore has order $n^2$. This order
is determined by the Alexander module and is unchanged by
replacing $t$ with $t^{-1}$, since
$f(t^{-1})=t^{-2}f(t)$. It is thus a knot invariant and
distinguishes the knots $J_n$ for $n>0$. In particular,
the subfamily $\{J_{2n}\}_{n\geq1}$ consists of pairwise
distinct knots.

Let $\widehat\varphi_n$ denote the monodromy of the capped
fiber of $J_n$. By the genus-two case of the Birman-Hilden
theorem \cite{BirmanHilden}, we may choose a representative commuting with a
fixed hyperelliptic involution; see
\cite[Section~9.4.2, p.~257]{FarbMargalit}.
The induced map on the quotient sphere defines a mapping class
\[
\overline\varphi_n\in\Mod(S^2,6),
\]
where the six marked points are the branch points of the
hyperelliptic covering.

\begin{lemma}\label{lem:three-three}
For every $n\geq0$, the mapping class $\overline\varphi_{2n}$
permutes the six branch points in two disjoint three-cycles.
\end{lemma}

\begin{proof}
We first consider the square knot. Its capped fiber is the
connected sum of the capped trefoil and mirror-trefoil fibers,
giving an orthogonal decomposition
\[
H_1(\Sigma_2;\mathbb F_2)=V_1\oplus V_2
\]
into two symplectic planes preserved by
$(\widehat\varphi_0)_*$. On each plane, the monodromy cyclically
permutes the three nonzero vectors, as follows from the
modulo-two trefoil calculation in Section~\ref{sec:four-strand};
the mirror monodromy has the same cycle type.

Recall that a quadratic refinement of the intersection pairing
is a function $q$ satisfying
\[
q(x+y)=q(x)+q(y)+x\cdot y.
\]
On a symplectic plane there is one odd refinement, which takes
the value one on every nonzero vector, and three even
refinements, each taking the value one on exactly one nonzero
vector. 
Consequently, for each $i=1,2$, the action induced by
$(\widehat\varphi_0)_*|_{V_i}$ on the quadratic refinements
of $V_i$ fixes the unique odd refinement and cyclically
permutes the three even refinements.

Every quadratic refinement $q$ on $V_1\oplus V_2$ is determined
by its restrictions to the two summands. Since the Arf invariant
is additive under orthogonal direct sum,
\[
\operatorname{Arf}(q)
=
\operatorname{Arf}(q|_{V_1})
+
\operatorname{Arf}(q|_{V_2})
\in\mathbb F_2.
\]
Thus $q$ is odd precisely when one restriction is odd and the
other is even. There are therefore six odd refinements,
forming two invariant triples: three with odd restriction
to $V_1$ and even restriction to $V_2$, and three with the
opposite parities. The action induced by
$(\widehat\varphi_0)_*$ cyclically permutes each triple.

We now identify this permutation with the permutation of the
branch points. Let
\[
B=\{b_1,\ldots,b_6\}
\]
be the branch set of the hyperelliptic covering.
The space $H_1(\Sigma_2;\mathbb F_2)$ can be identified with
the even-cardinality subsets of $B$, with each subset
identified with its complement. Under this identification,
addition is symmetric difference and the intersection
pairing is
\[
[S]\cdot[T]=|S\cap T|\pmod2.
\]
Every nonzero class is represented uniquely by a two-element
subset $\{b_j,b_k\}$; geometrically, it is the homology class
of the inverse image of an arc joining $b_j$ to $b_k$.
This description is compatible with the mapping class group
action; see \cite[Section~2]{McCarthyWeierstrass}.

For each $b_i\in B$, define
\[
q_i([S])
=
\frac{|S|}{2}+\epsilon_i(S)\pmod2,
\qquad
\epsilon_i(S)=
\begin{cases}
1,&b_i\in S,\\
0,&b_i\notin S.
\end{cases}
\]
Replacing $S$ by its complement gives
\[
\frac{|B\setminus S|}{2}+\epsilon_i(B\setminus S)
=
\frac{6-|S|}{2}+1-\epsilon_i(S)
\equiv
\frac{|S|}{2}+\epsilon_i(S)\pmod2,
\]
so $q_i$ is well defined. Moreover,
\[
|S\triangle T|=|S|+|T|-2|S\cap T|,
\qquad
\epsilon_i(S\triangle T)
\equiv\epsilon_i(S)+\epsilon_i(T)\pmod2.
\]
It follows that
\[
q_i([S]+[T])
=
q_i([S])+q_i([T])+[S]\cdot[T],
\]
so $q_i$ is a quadratic refinement. On the nonzero classes,
its values are
\[
q_i([\{b_j,b_k\}])
=
\begin{cases}
0,&i\in\{j,k\},\\
1,&i\notin\{j,k\}.
\end{cases}
\]

To check its parity, write $v_{jk}=[\{b_j,b_k\}]$.
The intersection formula shows that
$$(v_{12},v_{23},v_{45},v_{56})$$ is a symplectic basis for $H_1(\Sigma_2;\mathbb F_2)$.
Hence
\[
\operatorname{Arf}(q_1)
=
q_1(v_{12})q_1(v_{23})
+
q_1(v_{45})q_1(v_{56})
=
0\cdot1+1\cdot1=1.
\]
Relabeling the branch points proves that every $q_i$ is odd.
The six refinements are distinct because their sets of zeros
are distinct. Since there are exactly six odd refinements,
the assignment $b_i\mapsto q_i$ is a bijection.

Finally, let $h$ be a diffeomorphism commuting with the
hyperelliptic involution, and suppose its quotient sends
$b_i$ to $b_{\sigma(i)}$. The induced action on homology
sends $[S]$ to $[\sigma(S)]$. The defining formula therefore
gives
\[
q_{\sigma(i)}(h_*v)=q_i(v),
\mathrm{\ \ or\ equivalently\ \ }
q_{\sigma(i)}=q_i\circ h_*^{-1}.
\]
Thus the bijection between branch points and odd refinements
is equivariant. Applying this to $\widehat\varphi_0$, the
two three-cycles on the odd refinements give two three-cycles
on the branch points. Hence $\overline\varphi_0$ has the
required permutation.

For the even-indexed knots, use the monodromy in the proof of
\cite[Proposition~3.2, p.~490]{ParkYun}
with $(p,q)=(2n,0)$.
After capping, it gives
\[
\widehat\varphi_{2n}=T_c^{\,2n}\circ\widehat\varphi_0,
\]
where $T_c$ is the Dehn twist about the Stallings curve.
On modulo-two homology,
\[
(T_c)_*(v)=v+(v\cdot[c])[c].
\]
Since $[c]\cdot[c]=0$, applying this formula twice gives
$(T_c)_*^2=\mathrm{id}$. Thus
\[
(\widehat\varphi_{2n})_*=(\widehat\varphi_0)_*
\quad\text{on }H_1(\Sigma_2;\mathbb F_2).
\]
Their actions on the odd quadratic refinements, and hence
their permutations of the six branch points, agree.
\end{proof}

For each $n\geq0$, choose one of the two three-point orbits of
$\overline\varphi_{2n}$ and forget the branch points in that orbit.  
The induced mapping class in $\Mod(S^2,3)$ permutes the
remaining three points cyclically. Choose a disk-braid lift
$\gamma_n\in B_3$ of this mapping class and pull its configuration
loop back along the map $u:\Sigma_g\to S^1$ used in
Section~\ref{sec:braids}.
Proposition~\ref{prop:braid-surface} gives a connected symplectic
multisection $F_{\gamma_n}\subset X_g=\Sigma_g\times S^2$,
for every $g\geq1$, with
\[
[F_{\gamma_n}]=3[\Sigma_g\times\{\mathrm{pt}\}],
\qquad
[F_{\gamma_n}]^2=0,
\qquad
g(F_{\gamma_n})=3(g-1)+1.
\]

The distinction between the knots does not persist in the
three-point sphere mapping class. Indeed, the permutation map
\[
\Mod(S^2,3)\longrightarrow S_3
\]
is an isomorphism, and all three-cycles are conjugate.
Thus the reduced sphere mapping classes are conjugate for all
$n$. The knot distinction alone therefore does not establish
that the resulting degree-three multisections are distinct.
Any such conclusion would require additional information
about the chosen disk-braid lifts or the full six-point
monodromies.

\section{Braid variation in all five cases}
\label{sec:five-braids}

We now replace the reference multisections in the five pairs of
Section~\ref{sec:five-patterns} by braided multisections and vary
the surface gluing. For the product pairs, the replacements preserve
the cyclic covering maps, so the homology calculations of that
section continue to apply. The twisted degree-four multisection has
dihedral covering monodromy. We compute its covering image explicitly
and then prove that it too admits gluings with unbounded finite
first-homology orders.

\subsection{Braid replacements and the twisted multisections}

For a product summand with a degree-three multisection, choose a disk-braid lift of the three-point sphere mapping class obtained from the trefoil or
figure-eight hyperelliptic quotient, or from
Lemma~\ref{lem:three-three}. 
For a product summand with a degree-four multisection, use
\[
\beta_k=\sigma_1^{2k+1}\sigma_2\sigma_3,\qquad k\geq1.
\]
A product summand with a degree-two multisection may use either the
standard two-strand braid $\sigma_1$ or the trefoil braid $\sigma_1^3$.  

Recall that the cyclic covering $F_{g,p}\to\Sigma_g$ used in
Section~\ref{sec:five-patterns} corresponds to $\ker\chi_p$,
where $\chi_p(a_1)=1\in\Z/p$ and $\chi_p$ vanishes on the
other standard generators. For each of the above braid
multisections, number the sheets so that the monodromy
around $a_1$ is $(1\,2\,\cdots\,p)$. Since all other standard
base generators act trivially on the sheets, its permutation
monodromy agrees with that of $F_{g,p}\to\Sigma_g$.
The two coverings are therefore isomorphic over $\Sigma_g$,
and their pushforwards on first homology agree under the
induced identification of the covering surfaces.

We first construct the twisted degree-two multisection.
For $\theta\in\mathbb R$, let $R_x(\theta)$ and $R_z(\theta)$
denote rotations of the unit sphere $S^2\subset\mathbb R^3$
through angle $\theta$ about the $x$- and $z$-axes,
respectively. Consider the flat sphere bundle over $\Sigma_3$
with holonomies
\[
\rho(a_1)=R_x(\pi),\qquad \rho(b_1)=R_z(\pi),
\]
and trivial holonomy on the remaining standard generators.
Explicitly, these half-turns are
\[
R_x(\pi)(x,y,z)=(x,-y,-z),\qquad
R_z(\pi)(x,y,z)=(-x,-y,z).
\]
They commute in $SO(3)$, so the holonomies satisfy the
surface relation. Their lifts to $SU(2)$ have commutator
$-1$. 
Thus the associated oriented rank-three bundle has nonzero
second Stiefel-Whitney class, so its sphere bundle is the
nontrivial oriented $S^2$-bundle $E_3\to\Sigma_3$ introduced
in Section~\ref{sec:five-patterns}.

Let $N=(0,0,1)$ and $S=(0,0,-1)$ be the north and south
poles of $S^2$. In the flat-bundle description
\[
E_3=(\widetilde\Sigma_3\times S^2)/\pi_1(\Sigma_3),
\]
the action is
\[
\gamma\cdot(\widetilde y,z)
=
(\gamma\widetilde y,\rho(\gamma)z).
\]
Every holonomy preserves the unordered pair $\{N,S\}$,
so this action preserves
$\widetilde\Sigma_3\times\{N,S\}$. Its quotient defines
\[
A_{3,2}
=
(\widetilde\Sigma_3\times\{N,S\})/\pi_1(\Sigma_3)
\subset E_3.
\]
Over a sufficiently small disk $U\subset\Sigma_3$, a flat
trivialization identifies this subset with
\[
(U\times\{N\})\sqcup(U\times\{S\}).
\]
Thus $A_{3,2}$ is an embedded horizontal surface whose
projection to $\Sigma_3$ is an unbranched degree-two
covering. Since $\rho(a_1)=R_x(\pi)$ exchanges the poles,
the lift of $a_1$ joins the two points over the basepoint.
The covering is therefore connected.

The finite rotational holonomy induces finite holonomy
on the normal bundle of $A_{3,2}$, so this bundle becomes
trivial after pullback to a finite cover. Since pullback
multiplies the normal Euler number by the covering degree,
we obtain $[A_{3,2}]^2=0$. 
Because the multisection has degree two, we can write
$[A_{3,2}]=2[S_3]+b[D]$, where 
$S_3$ is a section and $D$\/ is a sphere fiber.  
Since $[S_3]^2=[S_3]\cdot[D]=1$ and $[D]^2=0$, we have 
\[
0=[A_{3,2}]^2=4+4b.
\]
Thus $b=-1$ and $[A_{3,2}]=2[S_3]-[D]$, while the unbranched covering formula gives
\[
g(A_{3,2})=2(3-1)+1=5.
\]

The holonomies preserve the standard area form on $S^2$.
Its sum with the pullback of a positive area form on
$\Sigma_3$ therefore descends from
$\widetilde\Sigma_3\times S^2$ to a symplectic form on $E_3$.
Since the sphere coordinate is locally constant along
$A_{3,2}$, the restriction of this form to $A_{3,2}$ is
the positive pullback of the base area form.
Consequently, $A_{3,2}$ is symplectic.

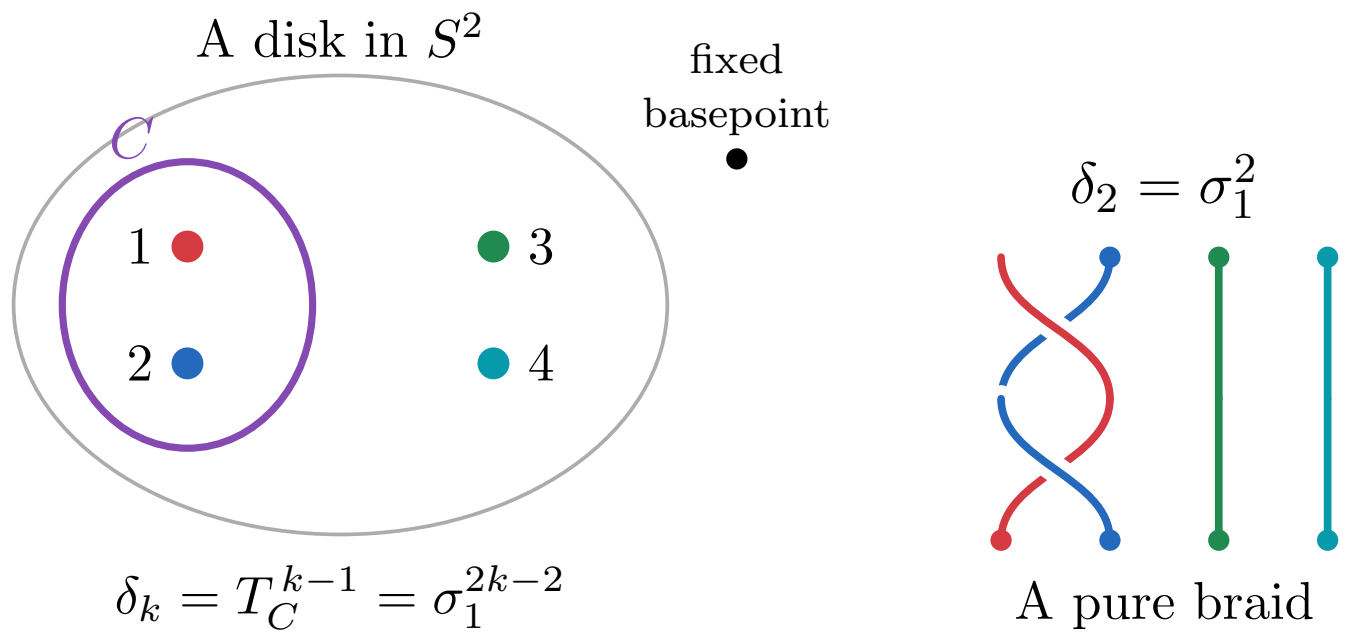
\begin{figure}[htbp]
\centering
\begin{tikzpicture}[x=0.92cm,y=0.92cm,line cap=round]
\draw[gray!65,line width=.7pt] (0,0) ellipse (2.35 and 1.65);
\node[font=\small] at (0,1.94) {A disk in $S^2$};
\draw[twistpurple,line width=1.35pt] (-1.1,0) ellipse (.9 and 1.03);
\node[text=twistpurple] at (-1.5,1.2) {$C$};
\fill[strandred] (-1.1,.42) circle (3pt);
\fill[strandblue] (-1.1,-.42) circle (3pt);
\fill[strandgreen] (1.1,.42) circle (3pt);
\fill[strandturquoise] (1.1,-.42) circle (3pt);
\node[left,font=\small] at (-1.2,.42) {$1$};
\node[left,font=\small] at (-1.2,-.42) {$2$};
\node[right,font=\small] at (1.2,.42) {$3$};
\node[right,font=\small] at (1.2,-.42) {$4$};
\fill[black] (2.85,1.05) circle (2pt);
\node[align=center,font=\scriptsize] at (2.85,1.55) {fixed\\basepoint};
\node[font=\small] at (0,-2.1) {$\delta_k=T_C^{\,k-1}=\sigma_1^{2k-2}$};
\end{tikzpicture}
\hspace{0.65cm}
\begin{tikzpicture}[x=0.72cm,y=0.72cm,line cap=round,line join=round]
\node at (1.5,0.65) {$\delta_2=\sigma_1^2$};
\fill[strandred] (0,0) circle (2pt);
\fill[strandblue] (1,0) circle (2pt);
\fill[strandgreen] (2,0) circle (2pt);
\fill[strandturquoise] (3,0) circle (2pt);
\draw[strandgreen,line width=1.45pt] (2,0.0000) -- (2,-1.3000);
\draw[strandturquoise,line width=1.45pt] (3,0.0000) -- (3,-1.3000);
\draw[strandblue,line width=1.45pt] (1,0.0000) .. controls (1,-0.5850) and (0,-0.7150) .. (0,-1.3000);
\draw[white,line width=5.1pt] (0,0.0000) .. controls (0,-0.5850) and (1,-0.7150) .. (1,-1.3000);
\draw[strandred,line width=1.45pt] (0,0.0000) .. controls (0,-0.5850) and (1,-0.7150) .. (1,-1.3000);
\draw[strandgreen,line width=1.45pt] (2,-1.3000) -- (2,-2.6000);
\draw[strandturquoise,line width=1.45pt] (3,-1.3000) -- (3,-2.6000);
\draw[strandred,line width=1.45pt] (1,-1.3000) .. controls (1,-1.8850) and (0,-2.0150) .. (0,-2.6000);
\draw[white,line width=5.1pt] (0,-1.3000) .. controls (0,-1.8850) and (1,-2.0150) .. (1,-2.6000);
\draw[strandblue,line width=1.45pt] (0,-1.3000) .. controls (0,-1.8850) and (1,-2.0150) .. (1,-2.6000);
\fill[strandred] (0,-2.6) circle (2pt);
\fill[strandblue] (1,-2.6) circle (2pt);
\fill[strandgreen] (2,-2.6) circle (2pt);
\fill[strandturquoise] (3,-2.6) circle (2pt);
\node[align=center,font=\small] at (1.5,-3.2) {A pure braid};
\end{tikzpicture}
\caption{The pure-braid insertion. The curve $C$ encloses the first two marked points; choose the sign of $T_C$ so that it represents $\sigma_1^2$. The right-hand panel shows one full twist. Its permutation is the identity, although its marked mapping class is nontrivial. The left-hand drawing is schematic; the four points may be moved to the equatorial configuration used in the flat-bundle construction.}
\label{fig:pure-insertion}
\end{figure}

Next we construct a family of degree-four multisections
in the nontrivial oriented sphere bundle $E_2\to\Sigma_2$.
Let
\[
\delta_k=\beta_k\beta_1^{-1}=\sigma_1^{2k-2},
\qquad k\geq1.
\]
These are pure braids, with $\delta_1=1$. 
Let $C$ be the boundary of a small regular neighborhood of the standard
arc joining the first two marked points and defining the half-twist
$\sigma_1$. Choose the sign of the Dehn twist $T_C$ so that it
represents $\sigma_1^2$.  
Thus $\delta_k=T_C^{k-1}$ in the marked sphere mapping class group;
see Figure~\ref{fig:pure-insertion}.

Let $\mathcal P_4$ be the four equally spaced equatorial points,
one of which is $(1,0,0)$. The reference flat model has holonomies
\[
\rho(a_1)=s=R_x(\pi),\quad
\rho(b_1)=c^2=R_z(\pi),\quad
\rho(a_2)=c=R_z(\pi/2),\quad
\rho(b_2)=1
\]
(cf.~\cite[Proposition~9.2]{AkhmedovMultisections}).
The standard cyclic braid
$\beta_0=\sigma_1\sigma_2\sigma_3$ has a spherical
representative given by a quarter-turn of four equally
spaced points. Choose an orientation-preserving
identification of this marked sphere with
$(S^2,\mathcal P_4)$ under which that representative
is $c=R_z(\pi/2)$. We use this identification for the
subsequent braid insertions.

Choose an area-preserving representative $d_k$ of $\delta_k$,
supported in an annulus about $C$ and equal to the identity near
$\mathcal P_4$, with $d_1=\mathrm{id}$. Replace the holonomies by
\[
\rho_k(a_1)=s,\qquad \rho_k(b_1)=c^2,\qquad
\rho_k(a_2)=d_kc,\qquad \rho_k(b_2)=1.
\]
They satisfy the surface relation because
$[s,c^2]=[d_kc,1]=1$.

After forgetting the marked points, the annular twist $d_k$ is
isotopic to the identity through area-preserving sphere
diffeomorphisms. 
Such an isotopy deforms $\rho_k$ to the reference
representation $\rho_1$, whose holonomies are the
rotations $s,c^2,c$, and the identity.
During the isotopy, only the $a_2$
holonomy changes, and the surface relation holds throughout.
The underlying sphere bundle is therefore unchanged. 
The rotations $s=R_x(\pi)$ and $c^2=R_z(\pi)$ are
half-turns about perpendicular axes, so any lifts of
them to $SU(2)$ have commutator $-1$.
Hence the sphere bundle associated to $\rho_k$ is $E_2\to \Sigma_2$ for every $k$.

The marked $a_2$ holonomy has spherical class
\[
[d_kc]=\delta_k\beta_0=\beta_{k-1}.
\]
For $k=1$, the $a_2$-holonomy is the original quarter-turn
$c=R_z(\pi/2)$. For $k=2$ and $k=3$, its spherical
mapping class is that of the stabilized trefoil and
cinquefoil braid, respectively.

\begin{proposition}\label{prop:twisted-braid}
For every $k\geq1$, the holonomy representation $\rho_k$
determines a connected degree-four multisection
$A_{2,4}^{(k)}\subset E_2$. This multisection is symplectic
and satisfies
\[
g(A_{2,4}^{(k)})=5,\qquad
[A_{2,4}^{(k)}]^2=0,\qquad
[A_{2,4}^{(k)}]=4[S_2]-2[D].
\]
The coverings $A_{2,4}^{(k)}\to\Sigma_2$ are mutually
isomorphic over $\Sigma_2$. With the identification of
the marked sphere fixed above, the mapping classes
$[\rho_k(a_2)]\in\Mod(S^2,4)$ are pairwise distinct,
although their permutations of the four marked points agree.
\end{proposition}

\begin{proof}
Using the diagonal action determined by $\rho_k$, define
\[
A_{2,4}^{(k)}
=
(\widetilde\Sigma_2\times\mathcal P_4)/\pi_1(\Sigma_2)
\subset E_2.
\]
Since $d_k$ fixes every point of $\mathcal P_4$, the
restrictions of $\rho_k$ and $\rho_1$ to this set agree.
In particular, the $a_2$-holonomy acts on $\mathcal P_4$
by the transitive four-cycle induced by $c$.

In a local flat trivialization over a disk
$U\subset\Sigma_2$, the subset $A_{2,4}^{(k)}$ is
$U\times\mathcal P_4$, a disjoint union of four horizontal
disks. Thus it is an embedded horizontal multisection,
and its projection is an unbranched degree-four covering.
The transitive action on $\mathcal P_4$ makes this
covering connected. Moreover, because the action on
$\mathcal P_4$ is independent of $k$, the identity on
$\widetilde\Sigma_2\times\mathcal P_4$ induces an
isomorphism
\[
A_{2,4}^{(k)}\longrightarrow A_{2,4}^{(1)}
\]
over $\Sigma_2$. The covering formula gives
\[
g(A_{2,4}^{(k)})=4(2-1)+1=5.
\]

We next compute the self-intersection. Since the multisection is horizontal, its normal bundle
is naturally identified with the vertical tangent bundle
restricted to it. Indeed, in a flat trivialization, each
sheet is locally $U\times\{p\}$, and its normal plane
is identified with $T_pS^2$. A sphere holonomy map $f$
induces normal parallel transport by
\[
Df_p\colon T_pS^2\longrightarrow T_{f(p)}S^2.
\]
Thus the normal holonomy is determined by the derivatives
of the sphere holonomies at the marked points.

By construction, $d_k$ is supported in an annular
neighborhood of $C$ disjoint from $\mathcal P_4$.
Since $c$ permutes $\mathcal P_4$, the map $d_k$
is the identity near $c(p)$ for every
$p\in\mathcal P_4$.
The chain rule therefore gives
\[
D(d_k\circ c)_p
=
(Dd_k)_{c(p)}\circ(Dc)_p
=
(Dc)_p.
\]
All other generator holonomies are unchanged.
Consequently, under the preceding identification of
the covering surfaces, the normal holonomy agrees
with that for $k=1$. The image of $\rho_1$ is the
finite rotation group generated by $s$ and $c$,
so this normal holonomy is finite.

The normal bundle therefore becomes trivial on a
finite cover of $A_{2,4}^{(k)}$. Since pullback
multiplies its Euler number by the covering degree,
we obtain $[A_{2,4}^{(k)}]^2=0$. The degree-four
condition gives
\[
[A_{2,4}^{(k)}]=4[S_2]+b[D].
\]
Using $[S_2]^2=[S_2]\cdot[D]=1$ and $[D]^2=0$, we find
\[
0=[A_{2,4}^{(k)}]^2=16+8b.
\]
Hence $b=-2$, proving the asserted homology class.

All holonomies preserve the standard sphere area form.
Its sum with the pullback of a positive area form on
$\Sigma_2$ descends to a symplectic form on $E_2$.
Along the horizontal multisection, the sphere area
form vanishes and the pullback of the base area form
is positive. Thus $A_{2,4}^{(k)}$ is symplectic.

Finally, in the fixed identification of the marked sphere,
\[
[\rho_k(a_2)]=[d_kc]=T_C^{\,k-1}[c].
\]
For any $k,\ell\geq1$, cancellation gives
\[
[\rho_k(a_2)]\,[\rho_\ell(a_2)]^{-1}
=
T_C^{\,k-\ell}.
\]
The curve $C$ separates the four marked points into
two pairs, so it is essential and $T_C$ has infinite
order in $\Mod(S^2,4)$. The displayed classes are
therefore distinct whenever $k\ne\ell$.
\end{proof}

\subsection{Homological data for the twisted and hybrid sums}

We next compute the image of the covering map
$A_{2,4}^{(k)}\to\Sigma_2$ on the first homology and
determine the meridian orders in the complements
of the twisted multisections. These calculations
provide the homological data needed for the
twisted and hybrid sums.

The permutations induced by $c$ and $s$ on
$\mathcal P_4$ generate the dihedral group
\[
D_4=\langle c,s\mid c^4=s^2=1,\ scs=c^{-1}\rangle
\]
of order eight. Let
$q:A_{2,4}^{(k)}\to\Sigma_2$ be the restriction
of the bundle projection. Since $d_k$ fixes
$\mathcal P_4$ pointwise, the permutation monodromy
of $q$ is
\[
\bar\rho:\pi_1(\Sigma_2)\longrightarrow D_4,
\qquad
a_1\longmapsto s,\quad b_1\longmapsto c^2,\quad
a_2\longmapsto c,\quad b_2\longmapsto1.
\]

Choose a basepoint in $\Sigma_2$ and identify its
fiber with $\mathcal P_4$. A based loop lifts to a
closed loop starting at the point corresponding
to $p_0=(1,0,0)$ precisely when its permutation
monodromy fixes $p_0$. Since the stabilizer of
$p_0$ in $D_4$ is $\langle s\rangle$, the covering
subgroup is
\[
\bar\rho^{-1}(\langle s\rangle)
\subset\pi_1(\Sigma_2).
\]
After identifying $A_{2,4}^{(k)}$ with $\Sigma_5$,
write
\[
Q^{D_4}_{4,2}=q_*:
H_1(\Sigma_5;\Z)\longrightarrow H_1(\Sigma_2;\Z).
\]
We use $a_i,b_i$ also to denote the corresponding
homology classes in $\Sigma_2$.

Consider the homomorphism $D_4\to\Z/2$ sending
$c$ to $1$ and $s$ to $0$. Its composite with
$\bar\rho$ is reduction modulo two of the
$a_2$-coordinate. Since this composite vanishes
on the covering subgroup,
$\operatorname{im}Q^{D_4}_{4,2}$ is contained in
the sublattice with even $a_2$-coordinate.
For the reverse inclusion, observe that the
covering subgroup contains
\[
a_1,\qquad b_2,\qquad
b_1[a_1,a_2]^{-1},\qquad a_2^2b_1^{-1}.
\]
The first two elements have permutation images
$s$ and $1$. Moreover, $[s,c]=c^2$, so the last
two elements have trivial permutation image.
Their homology classes are
\[
a_1,\qquad b_2,\qquad b_1,\qquad 2a_2-b_1,
\]
which generate the entire sublattice with even
$a_2$-coordinate. Therefore
\[
\operatorname{im}Q^{D_4}_{4,2}
=\langle a_1,b_1,2a_2,b_2\rangle,
\qquad
\operatorname{coker}Q^{D_4}_{4,2}\cong\Z/2.
\]

For the cyclic degree-four covering,
$\operatorname{coker}Q_{4,2}\cong\Z/4$.
Thus the two covering pushforwards cannot be
identified by integral changes of bases.
Consequently, the cyclic degree-four calculations
used in \cite[Theorem~8.2 and Proposition~9.4]{AY}
cannot be transferred to the present covering
by such changes of bases. By
Proposition~\ref{prop:twisted-braid}, the
pushforwards $Q^{D_4}_{4,2}$ agree for all $k$
under the induced identifications of the
covering surfaces.

For the complement calculation, let
$A_{g,p}\subset E_g$ be one of the twisted
multisections and put
\[
M=E_g\setminus\Int\nu A_{g,p}.
\]
The exterior $M$ is a $P_p$-bundle over
$\Sigma_g$. Choose lifts $A_i,B_i\in\pi_1(M)$
of the standard base generators and based
puncture loops $x_1,\ldots,x_p$ in the fiber.
Let $\Phi_i$ and $\Psi_i$ be the automorphisms
of $\pi_1(P_p)$ induced by conjugation by
$A_i$ and $B_i$, respectively.
The product $\prod_{i=1}^g[A_i,B_i]$ projects
to the identity in $\pi_1(\Sigma_g)$, so it
equals a word $w\in\pi_1(P_p)$. Hence
$\pi_1(M)$ has a presentation of the form
\[
\left\langle A_i,B_i,x_j\ \middle|\
\begin{gathered}
A_ix_jA_i^{-1}=\Phi_i(x_j),\quad
B_ix_jB_i^{-1}=\Psi_i(x_j),\\
x_1\cdots x_p=1,\quad \prod_{i=1}^g[A_i,B_i]=w
\end{gathered}\;\right\rangle,
\]
where $1\leq i\leq g$ and $1\leq j\leq p$.
Conjugating the fiber subgroup by the last
relation gives
\[
\prod_{i=1}^g[\Phi_i,\Psi_i]
=\operatorname{Inn}(w),
\qquad
\operatorname{Inn}(w)(u)=wuw^{-1}.
\]
The automorphisms and the word $w$ depend on
the bundle and the chosen lifts.

Although we do not need the explicit word $w$,
we can determine its contribution to the
abelianization. The multisection class is
\[
[A_{g,p}]=p[S_g]-(p/2)[D],
\]
and its intersections with $S_g$ and $D$ are
$p/2$ and $p$, respectively. The homology
sequence of the pair, as in
Section~\ref{sec:five-patterns}, therefore
shows that the meridian has order $p/2$.
Since $H_1(\Sigma_g;\Z)$ is free abelian,
the resulting sequence splits, giving a
noncanonical isomorphism
\[
H_1(M;\Z)
\cong H_1(\Sigma_g;\Z)\oplus\Z/(p/2),
\]
with the meridian generating the cyclic summand.

We now compare this with the abelianization
of the displayed presentation. The permutation
monodromy is transitive, so the conjugation
relations identify all puncture classes with
a single meridian class $\mu$. The relation
$x_1\cdots x_p=1$ then gives $p\mu=0$.
Thus, before imposing the surface relation,
the puncture classes generate $\Z/p$.
Let $[w]$ denote the image of $w$ in this
cyclic group. The surface relation imposes
$[w]=0$, so the resulting meridian subgroup is
\[
(\Z/p)/\langle[w]\rangle.
\]
Its order is $p/2$, and hence
$\langle[w]\rangle$ has order two.
Since $\Z/p$ has a unique element of order two,
we obtain
\[
[w]=(p/2)\mu\quad\text{in }\Z/p.
\]
In particular, $[w]=\mu$ when $p=2$ and
$[w]=2\mu$ when $p=4$. Setting $w=1$ would
therefore give the wrong abelianization in
both cases. 
This calculation determines the class $[w]\in\Z/p$,
but does not determine the word $w$ or the full
fundamental-group extension.

\subsection{Variation of the gluing}

We use the following linear-algebra argument to treat the actual
covering maps, without requiring cyclic monodromy.

\begin{proposition}\label{prop:unbounded-gluing}
Let $q_i:\Sigma_G\to\Sigma_{g_i}$ be connected unbranched
covering maps, with $G>1$ and $g_1+g_2=G$. Write $Q_i=(q_i)_*$
in symplectic integral bases. Then
\[
D(S)=\det\begin{pmatrix}Q_1\\-Q_2S\end{pmatrix},
\qquad S\in\operatorname{Sp}(2G,\Z),
\]
has unbounded absolute value. 
In particular, one can choose a sequence
$S_m\in\operatorname{Sp}(2G,\Z)$, $m\geq1$, such that the block matrices
\[
\begin{pmatrix}
Q_1\\
-Q_2S_m
\end{pmatrix}
\]
are nonsingular and their cokernels have strictly increasing finite orders.
\end{proposition}

\begin{proof}
We work first in $V=H_1(\Sigma_G;\Q)$, equipped with
its intersection pairing. Let
\[
\operatorname{PD}_G:
H^1(\Sigma_G;\Q)\longrightarrow H_1(\Sigma_G;\Q)
\]
be the Poincar\'e duality isomorphism, and set
\[
W_i=\operatorname{PD}_G(\operatorname{im}q_i^*).
\]
For $\alpha,\beta\in H^1(\Sigma_{g_i};\Q)$, we have
\[
\bigl\langle q_i^*\alpha\cup q_i^*\beta,
[\Sigma_G]\bigr\rangle
=
\deg(q_i)\bigl\langle\alpha\cup\beta,
[\Sigma_{g_i}]\bigr\rangle.
\]
Since $\deg(q_i)\ne0$ and the cup-product pairing
on $H^1(\Sigma_{g_i};\Q)$ is nondegenerate,
the last identity implies that $q_i^*$ is
injective. It also shows that $W_i$ is a symplectic
subspace of dimension $2g_i$.

For $x\in V$ and $\alpha\in H^1(\Sigma_{g_i};\Q)$,
Poincar\'e duality and the compatibility of pullback
and pushforward with evaluation give
\[
\bigl\langle x,\operatorname{PD}_G(q_i^*\alpha)\bigr\rangle
=(q_i^*\alpha)(x)
=\alpha(Q_i x).
\]
Consequently, $x$ is orthogonal to $W_i$ precisely
when $\alpha(Q_i x)=0$ for every $\alpha$, or
equivalently when $Q_i x=0$. Hence
\[
K_i=\ker Q_i=W_i^\perp.
\]
The symplectic orthogonal complement of a symplectic
subspace is symplectic, so each $K_i$ is symplectic
and has dimension $2G-2g_i$. Since $g_1+g_2=G$,
\[
\dim K_1=2g_2,\qquad \dim K_2=2g_1.
\]

Moreover, the square matrix defining $D(S)$ has kernel
\[
\ker\begin{pmatrix}Q_1\\-Q_2S\end{pmatrix}
=K_1\cap S^{-1}(K_2).
\]
Applying $S$ to this kernel gives
$S(K_1)\cap K_2$. Hence
\[
D(S)\ne0
\quad\Longleftrightarrow\quad
S(K_1)\cap K_2=\{0\}.
\]

The symplectic subspaces $K_1$ and $K_2^\perp$
both have dimension $2g_2$. Choose rational
symplectic bases for these subspaces and extend
them to symplectic bases for $V$. The linear map
sending the first extended basis to the second
gives $S_+\in\operatorname{Sp}(2G,\Q)$ with
\[
S_+(K_1)=K_2^\perp.
\]
Since $K_2$ is symplectic,
$K_2^\perp\cap K_2=\{0\}$, and therefore
$D(S_+)\ne0$.

On the other hand, choose nonzero vectors
$x\in K_1$ and $y\in K_2$. Extending each vector
to a rational symplectic basis of $V$ gives
$S_0\in\operatorname{Sp}(2G,\Q)$ with $S_0x=y$.
Thus $y$ is a nonzero element of
$S_0(K_1)\cap K_2$, so $D(S_0)=0$.
Taking $R=S_+$ if $D(I)=0$, and $R=S_0$
otherwise, we obtain
\[
R\in\operatorname{Sp}(2G,\Q),
\qquad D(R)\ne D(I).
\]

By \cite[Theorem~8.5]{Taylor}, the group
$\operatorname{Sp}(V)$ is generated by the
symplectic transvections
\[
\tau_v(t)(x)=x+t\langle x,v\rangle v,
\qquad v\in V,\quad t\in\Q.
\]
Consequently, after rescaling each direction to be integral, write
\[
R=\tau_{v_1}(r_1)\cdots\tau_{v_\ell}(r_\ell),
\qquad v_j\in H_1(\Sigma_G;\Z),\quad r_j\in\Q.
\]
Choose a positive integer $N$ such that every $Nr_j$ is integral,
and set
\[
S(t)=\tau_{v_1}(Nr_1t)\cdots\tau_{v_\ell}(Nr_\ell t).
\]
Its entries are integer polynomials in $t$, and $S(n)$ belongs to
$\operatorname{Sp}(2G,\Z)$ for every integer $n$. The polynomial
$D(S(t))$ is nonconstant because
\[
D(S(0))=D(I)\ne D(R)=D(S(1/N)).
\]
Its absolute value is therefore unbounded on the positive integers.
Choose a subsequence of nonzero values with strictly increasing
absolute values. The cokernel of a nonsingular square integral
matrix has order equal to the absolute value of its determinant,
which proves the assertion.
\end{proof}

After matching the areas of the summing surfaces,
choose an area-preserving gluing map $\psi_S$
from the first summing surface to the second
whose induced map on the first homology is $S$ in
the chosen bases, as in
Section~\ref{sec:five-patterns}.
Choose adapted normal framings as in
Lemma~\ref{lem:peripheral}, and let $Z=Z_S$ denote
the symplectic sum obtained using $\psi_S$ with
zero normal twist relative to these framings.
Formula~\eqref{eq:homology-formula} then gives
\[
H_1(Z_S;\Z)\cong
\operatorname{coker}\begin{pmatrix}Q_1\\-Q_2S\end{pmatrix}
\oplus\Z/\gcd(d_1,d_2),
\]
where $d_i$ is the meridian order in the $i$th exterior.
For the twisted case $(d_1,d_2)=(1,2)$; for the hybrid case it is
$(2,2)$. Thus, when the determinant is nonzero, the respective
first-homology orders are $|D(S)|$ and $2|D(S)|$, using
$Q_2=Q^{D_4}_{4,2}$ in both cases. Proposition~\ref{prop:unbounded-gluing}
therefore supplies the required sequences for these two cases.
Because the covering maps and meridian orders are independent of
$k$, the sequence of homological gluings can be chosen once for
all braid parameters in each case.

\begin{theorem}\label{thm:five-families}
For each of the five pairs in
Section~{\rm \ref{sec:five-patterns}}, fix braid inputs
as described above. Each pair admits an infinite
family of closed, aspherical, minimal symplectic
sums $\{Z_m\}_{m\geq1}$,
obtained by varying the gluing while keeping
the braid inputs fixed. 
Every member satisfies
\[
e(Z_m)=4,\qquad \sigma(Z_m)=0,\qquad b_1(Z_m)=0.
\]
Within each family, the first homology groups
are finite and have strictly increasing orders.
Consequently, the manifolds in each family
are pairwise nondiffeomorphic.

For the three product pairs, the explicit
gluings of Section~{\rm \ref{sec:five-patterns}} give
\[
\begin{array}{c|c}
\text{pair}&H_1(Z_m;\Z)\\ \hline
\text{product }(2,3)&\Z/(6m)\\
\text{product }(2,4)&(\Z/2)^2\oplus\Z/(4m)\\
\text{product }(3,3)&(\Z/3)^2\oplus\Z/(3m).
\end{array}
\]
In each of the three $(2,4)$ cases, the
four-strand braid parameter gives a family
$Z_{k,m}$, where $k,m\geq1$. In the product $(2,4)$ case,
$k$ specifies the multisection
$F_{\beta_k}\subset X_2$, where
\[
\beta_k=\sigma_1^{2k+1}\sigma_2\sigma_3.
\]
In the twisted and hybrid $(2,4)$ cases, $k$\/ specifies
the multisection $A_{2,4}^{(k)}\subset E_2$
of Proposition~{\rm \ref{prop:twisted-braid}}.
In each $(2,4)$ case, we choose the degree-two
multisection once and for all. The parameter $k$
varies only the degree-four multisection,
while $m$ indexes the chosen gluings.
Therefore, for every fixed $k$, each $(2,4)$ case has its own family
$\{Z_{k,m}\}_{m\geq1}$ satisfying the conclusions
above, and the isomorphism type of $H_1(Z_{k,m};\Z)$ depends
only on $m$.
\end{theorem}

\begin{proof}
The product braid replacements preserve the
covering pushforwards and the multisection classes.
Lemma~\ref{lem:peripheral} provides adapted normal
framings for which the induced boundary maps on
the first homology agree with those of the original
multisections. Thus the gluing matrices and
homology calculations in
Section~\ref{sec:five-patterns} give the three
displayed groups. 

For the twisted and hybrid pairs,
Proposition~\ref{prop:twisted-braid} identifies
the degree-four coverings for all $k$ and shows
that their multisection classes agree.
Consequently, the covering pushforwards and
meridian orders are independent of $k$.
Since $g_1+g_2=G$ in both cases,
Proposition~\ref{prop:unbounded-gluing} and
formula~\eqref{eq:homology-formula} provide
gluing matrices for which the first homology
groups have strictly increasing finite orders.
These matrices can be chosen independently
of $k$. The same independence holds for the
explicit matrices in the product $(2,4)$ case.

As explained above, each chosen gluing matrix
is realized by an area-preserving diffeomorphism
after matching the areas of the summing surfaces.
Using the adapted normal framings and zero
normal twist therefore gives symplectic sums.
The sum formulas give $e=4$ and $\sigma=0$,
while finiteness of first homology gives $b_1=0$.
Asphericity and minimality follow from
Proposition~\ref{prop:aspherical-sums}.
Finally, the strictly increasing first-homology
orders distinguish the manifolds in each sequence.
\end{proof}

\begin{corollary}\label{cor:nonabelian}
For every product sum using the braid complements above, both
exterior groups inject into the closed sum group, and
\[
\pi_1(\Sigma_G)\times\Z\ \hookrightarrow\ \pi_1(Z).
\]
In each of the three product pairs of
Theorem~{\rm \ref{thm:five-families}}, the fundamental groups are pairwise
nonisomorphic as $m$ varies with the braid fixed.
\end{corollary}

\begin{proof}
The boundary maps are injective by
Proposition~\ref{prop:peripheral-presentation}. The normal-form
theorem for amalgamated products therefore embeds the exterior
groups and the common boundary group
$\pi_1(\Sigma_G\times S^1)$ in $\pi_1(Z)$.
The abelianizations have orders $6m$, $16m$, and $27m$ in the
three product rows, so they distinguish the groups within each row.
\end{proof}

\subsection{Noncomplex families of symplectic sums}

\begin{theorem}\label{thm:noncomplex}
Every manifold constructed in Theorem~{\rm \ref{thm:five-families}} has the rational
cohomology ring of $S^2\times S^2$.
For each of the five pairs, all but finitely
many members of the family $\{Z_m\}_{m\geq1}$
admit no complex structure.
Moreover, in each of the three $(2,4)$ cases,
there is an integer $m_0$, independent of $k$,
such that $Z_{k,m}$ admits no complex structure
for every $k\geq1$ and $m\geq m_0$.
Consequently, each pair yields infinitely many
pairwise nondiffeomorphic, aspherical, minimal
symplectic four-manifolds with the rational
cohomology ring of $S^2\times S^2$ and admitting
no complex structure.
\end{theorem}

\begin{proof}
Let $Z$ be any manifold constructed in
Theorem~\ref{thm:five-families}.
Since $b_1(Z)=0$ and $e(Z)=4$, Poincar\'e
duality gives $b_2(Z)=2$. Its signature is zero,
so its intersection form on
$H_2(Z;\Z)/\mathrm{Tor}$ is indefinite,
unimodular, and of rank two. It is therefore
equivalent over $\Z$ to one of
\[
\begin{pmatrix}0&1\\1&0\end{pmatrix},
\qquad
\begin{pmatrix}1&0\\0&-1\end{pmatrix}.
\]
Both forms are equivalent over $\Q$ to the
intersection form of $S^2\times S^2$.
Since $H^1(Z;\Q)=H^3(Z;\Q)=0$, it follows that
\[
H^*(Z;\Q)\cong H^*(S^2\times S^2;\Q)
\]
as graded rings.

Suppose now that $Z$ admits a complex structure.
Since $b_1(Z)=0$, the resulting complex surface
is K\"ahler by \cite[Theorem~11]{Buchdahl}.
Asphericity excludes exceptional spheres and rational or
ruled surfaces. Thus this complex surface is
minimal and has nonnegative Kodaira dimension.

Use the orientation induced by the complex
structure. The Euler characteristic remains
$4$, and the signature is still zero, regardless
of whether this orientation agrees with the
symplectic orientation. The Chern-number
identities therefore give
\[
\chi(\mathcal O_Z)
=\frac{e(Z)+\sigma(Z)}4=1,
\qquad
K_Z^2=2e(Z)+3\sigma(Z)=8.
\]
A minimal complex surface of Kodaira dimension
zero or one has $K^2=0$. Hence $Z$ must be a
minimal surface of general type.

Minimal surfaces of general type with fixed
invariants $(K^2,\chi(\mathcal O))=(8,1)$ have
only finitely many differentiable types
\cite[Corollary~165 on p.~402]{CataneseModuli}.
Consequently, there is a constant $B$ bounding
the orders of their first homology groups
whenever these groups are finite.
Thus any manifold in our families satisfying
\[
|H_1(Z;\Z)|>B
\]
admits no complex structure.

In each family $\{Z_m\}$, the first-homology
orders increase strictly and therefore
eventually exceed $B$. This proves the
noncomplexity assertion for all sufficiently
large $m$. In each $(2,4)$ case,
Theorem~\ref{thm:five-families} also shows that
$H_1(Z_{k,m};\Z)$ is independent of $k$.
The same cutoff $m_0$ therefore works for
every $k\geq1$ within that case.
The remaining conclusions follow from
Theorem~\ref{thm:five-families}.
\end{proof}

\begin{remark}
In each of the product, twisted, and hybrid $(2,4)$
families, the sequence $\{Z_{m,m}\}_{m\geq1}$ varies
both the braid and gluing parameters. Its members
have strictly increasing first-homology orders and
are therefore pairwise nondiffeomorphic. By
Theorem~\ref{thm:noncomplex}, all sufficiently large
members admit no complex structure.
\end{remark}

\begin{remark}\label{rem:spin}
Every twisted or hybrid sum is nonspin, independently
of the braid parameter and the gluing. Indeed, set
\[
M=E_2\setminus\Int\nu A_{2,4}^{(k)}.
\]
Since
\[
[A_{2,4}^{(k)}]=4[S_2]-2[D],
\]
the mod-two homology class of the multisection
vanishes. 
Applying the mod-two cohomological version of
the Thom-isomorphism argument used in
Section~\ref{sec:five-patterns}, the exact sequence
of the pair $(E_2,M)$ shows that the kernel of
\[
H^2(E_2;\Z/2)\longrightarrow H^2(M;\Z/2)
\]
is generated by the mod-two Poincar\'e dual of
$[A_{2,4}^{(k)}]$. The restriction map is therefore
injective.

On the other hand,
\[
\bigl\langle w_2(E_2),[S_2]\bigr\rangle
=S_2^2=1\pmod 2,
\]
so $w_2(E_2)$ restricts nontrivially to $M$.
For any sum $Z$ containing this complement,
naturality gives
\[
w_2(Z)|_M=w_2(M)=w_2(E_2)|_M\ne0.
\]
Thus $Z$ is nonspin.

Although these sums are nonspin, this alone does
not show that their free integral intersection
forms are odd: when $H_1(Z;\Z)$ has two-torsion,
an even intersection form need not imply that
$Z$ is spin.
For a product sum, both complements are spin,
but the sum is spin precisely when spin structures
on the two complements can be chosen so that
their boundary restrictions agree under the
gluing map. The preceding homology calculations
do not establish this compatibility.
\end{remark}

\section{Distinguishing symplectic sums with the same gluing}
\label{sec:fixed-gluing}

Throughout this section we consider the $(p,q)=(2,4)$ product pair
of Section~\ref{sec:five-patterns}.  
For $k\geq1$, consider the sum of $\Sigma_2\times S^2$ and
$\Sigma_3\times S^2$ along the genus-five multisections
$F_{\beta_k}$ and $F_{3,2}$ of Proposition~\ref{prop:four-sum}.
The multisection $F_{\beta_k}$ is defined by the four-strand braid
$\beta_k=\sigma_1^{2k+1}\sigma_2\sigma_3$, and $F_{3,2}$ is defined
by the positive half-twist $\sigma_1\in B_2$.
To compare the sums for different $k$, we first
specify the surface generators and their lifts to the two boundaries.

Example~\ref{ex:beta_k} gives presentations of the two complement
groups $G_{\beta_k,2}$ and $G_{\sigma_1,3}$. On the four-sheeted side
we use the generators $x_1,\ldots,x_4,a,b,c,d$ of that example.
On the two-sheeted side we denote the three pairs of base generators
by $(\alpha,\eta),(\gamma,\delta),(\epsilon,\zeta)$.
The chosen meridians are $x_1$ and $t$, respectively.

The base loop $a^4$ lifts to a closed loop on $F_{\beta_k}$ because
$a$ cyclically permutes the four sheets. To write the gluing relation
for this surface generator, we express a push-off of the lifted loop
in the generators of the exterior group.
Write $\Phi=\Phi_{\beta_k}$,
$\Delta=x_1x_2$, and $z=\Delta^kx_1\Delta^{-k}$. In the Artin
formulas of Example~\ref{ex:beta_k}, $\Phi$ takes each meridian
to a conjugate of the next in the cycle $(1\,2\,3\,4)$.
The conjugating words may be taken, in order, to be
$\Delta^kx_1,z,z,\Delta^kx_1$. Iterating gives
\[
 \Phi^4(x_1)=U_kx_1U_k^{-1},\quad \mathrm{where\ \ }
 U_k=\Phi^3(\Delta^kx_1)\,\Phi^2(z)\,\Phi(z)\,\Delta^kx_1.
\]
Since $a^4x_1a^{-4}=\Phi^4(x_1)=U_kx_1U_k^{-1}$, conjugating both sides by $U_k^{-1}$ gives $Ax_1A^{-1}=x_1$, where $A=U_k^{-1}a^4$.
Thus $A=U_k^{-1}a^4$ centralizes $x_1$. Before imposing the
sphere relation $x_1x_2x_3x_4=1$, the total exponent sum of $U_k$
is $(2k+1)+1+1+(2k+1)=4k+4$, since $\Phi$ preserves total
exponent sum. Hence $A$ has zero meridian coordinate in the
$\Z/4$ summand of Lemma~\ref{lem:h1}.
For the stabilized trefoil, substituting $x_4=(x_1x_2x_3)^{-1}$
and freely reducing gives
\[
 U_1=x_1x_2x_3x_1x_2^{-1}x_1^{-1}x_3^{-1}x_1^{-1}x_3^{-1}x_1.
\]

On the four-sheeted side set $B=b$ and
\[
 C_r=a^rca^{-r},\qquad D_r=a^rda^{-r}\qquad(0\leq r<4).
\]
On the two-sheeted side set $A'=\alpha^2$, $B'=\eta$, and, for
$r=0,1$, set
\[
\begin{aligned}
 C'_r&=\alpha^r\gamma\alpha^{-r},&
 D'_r&=\alpha^r\delta\alpha^{-r},\\
 E'_r&=\alpha^r\epsilon\alpha^{-r},&
 H'_r&=\alpha^r\zeta\alpha^{-r}.
\end{aligned}
\]
With the distinguished handle first, the ordered tuples of lifts are
\[
\begin{aligned}
 \mathcal L_k&=(A,B,C_0,D_0,C_1,D_1,C_2,D_2,C_3,D_3),\\
 \mathcal L'&=(A',B',C'_0,D'_0,E'_0,H'_0,C'_1,D'_1,E'_1,H'_1).
\end{aligned}
\]
Let $s_j\in\pi_1(F_{\beta_k})$ and $s'_j\in\pi_1(F_{3,2})$ be the
images of $(\mathcal L_k)_j$ and $(\mathcal L')_j$, respectively,
under the boundary projections. These are the standard generators
of the two genus-five surface groups. On each side, the chosen lifts
satisfy the surface relation, as verified below. Sending each
generator to its lift therefore defines a homomorphism from the
surface group to the boundary group whose composite with the
projection is the identity. These homomorphisms determine our
normal framings.

Let $\psi:F_{\beta_k}\to F_{3,2}$ be a based orientation-preserving
diffeomorphism. For $1\leq j\leq10$, write
$\psi_*(s_j)=w_j(s'_1,\ldots,s'_{10})$. Substituting the chosen
lift $(\mathcal L')_i$ for each $s'_i$ in $w_j$ gives the peripheral
lift $w_j(\mathcal L')$ of $\psi_*(s_j)$. Denote by $Z_{k,\psi}$
the sum formed with surface gluing $\psi$ and zero normal twist
in these framings.

\begin{proposition}\label{prop:explicit-closed-group}
The group $\pi_1(Z_{k,\psi})$ is obtained from the free product
$G_{\beta_k,2}*G_{\sigma_1,3}$ by imposing the relations
\[
 x_1t=1,\qquad (\mathcal L_k)_j=w_j(\mathcal L')
 \quad(1\leq j\leq10).
\]
\end{proposition}

\begin{proof}
The four-sheeted covering identifies $\pi_1(F_{\beta_k})$ with
$\ker\chi_4\subset\pi_1(\Sigma_2)$. Using the coset representatives
$1,a,a^2,a^3$, the nontrivial Schreier generators are $a^4$ and
$b_r=a^rba^{-r},C_r,D_r$ for $0\leq r<4$. Set $b_4=a^4ba^{-4}$.
Conjugating the base relator $[a,b][c,d]=1$ by $a^r$ gives
\[
 b_{r+1}=[C_r,D_r]^{-1}b_r\qquad(0\leq r<4).
\]
The first three equations eliminate $b_1,b_2,b_3$. Since $b_0=b$,
the fourth leaves the single relation
\[
 [a^4,b]\prod_{r=0}^3[C_r,D_r]=1.
\]

We now replace $a^4$ by its peripheral lift $A=U_k^{-1}a^4$.
We use $a,b,c,d$ to denote their lifts to $G_{\beta_k,2}$
through the constant section $\Sigma_2\times\{z_0\}$.
By Proposition~\ref{prop:group-presentation}, the lifts of $b,c,d$
commute with each of the fiber generators $x_1,\ldots,x_4$ and
hence with the entire fiber group.
For any fiber element $u$, normality of the fiber group implies
that $a^{-r}ua^r$ is again a fiber element. Thus $b,c,d$ commute
with $a^{-r}ua^r$. Conjugating by $a^r$ shows that $b_r,C_r,D_r$
commute with $u$.
In particular, $b_4$ commutes with $U_k$, and hence
\[
 ABA^{-1}=U_k^{-1}b_4U_k=b_4,
 \qquad [A,B]=[a^4,b].
\]
Thus the entries of $\mathcal L_k$ satisfy the same surface relation.
They also centralize $x_1$, by the choice of $A$ and the observation
above for $B,C_r,D_r$.

For the two-sheeted cover, let $R_r=[C'_r,D'_r][E'_r,H'_r]$ for
$r=0,1$, so that $R_1=\alpha R_0\alpha^{-1}$. The base relation
$[\alpha,\eta]R_0=1$ gives $\alpha\eta\alpha^{-1}=R_0^{-1}\eta$.
Conjugating by $\alpha$ and substituting gives
\[
 \alpha^2\eta\alpha^{-2}
 =R_1^{-1}(\alpha\eta\alpha^{-1})
 =R_1^{-1}R_0^{-1}\eta.
\]
Thus eliminating $\alpha\eta\alpha^{-1}$ from the two lifted
relations leaves $[A',B']R_0R_1=1$, since $A'=\alpha^2$ and
$B'=\eta$. In the exterior group, the lifts of
$\eta,\gamma,\delta,\epsilon,\zeta$ commute with $t$.
Since $\alpha$ inverts $t$, its square centralizes $t$ and
conjugation by $\alpha$ preserves the centralizer of $t$.
Hence every entry of $\mathcal L'$ centralizes $t$.

By Proposition~\ref{prop:peripheral-presentation}, the centralizer
of each meridian is the corresponding boundary group. 
On each side, the chosen lifts satisfy the defining relation of the
genus-five surface group. Thus the assignment of each surface
generator to its lift extends to a homomorphism from the surface
group to the corresponding boundary group.
This choice determines a normal framing
up to homotopy. Indeed, relative to any initial framing, the lifts
differ by meridian powers given by a homomorphism
$\pi_1(\Sigma_5)\to\Z$, which is realized by a map $\Sigma_5\to S^1$
rotating that framing.

In these framings the boundary gluing is $(s,z)\mapsto(\psi(s),\bar z)$.
It identifies $x_1$ with $t^{-1}$ and $(\mathcal L_k)_j$ with
$w_j(\mathcal L')$ for $1\leq j\leq10$. The meridian and the ten
surface lifts generate the boundary group, so the Seifert-Van Kampen
theorem gives the stated presentation.
\end{proof}

More generally, a normal twist $\theta:\Sigma_5\to S^1$ replaces
the $j$-th surface relation by
$(\mathcal L_k)_j=w_j(\mathcal L')t^{f_j}$, where $f_j$ is the
value of $\theta_*:H_1(\Sigma_5;\Z)\to\Z$ on the $j$-th generator.
For the rest of the section we take $\theta\equiv 1$. For each fixed $m$,
we choose the based gluings $\psi=\psi_m:F_{\beta_k}\to F_{3,2}$ so that
the words representing $(\psi_m)_*(s_1),\ldots,(\psi_m)_*(s_{10})$
in the target generators $s'_1,\ldots,s'_{10}$ are independent of $k$.

Write 
\((u_i,v_i)\), \(0\leq i\leq4\), for the five successive pairs of
\(\mathcal L'\). 
For each $c\in\mathbb Z$ and $0\leq i\leq4$, define the following
substitutions on the pair $(u_i,v_i)$:
\[
\mathcal U_c(u_i)=u_i,\quad \mathcal U_c(v_i)=v_i u_i^c,
\qquad
\mathcal L_c(u_i)=u_i v_i^c,\quad \mathcal L_c(v_i)=v_i.
\]
Each substitution fixes the generators in the other four pairs, 
and preserves \([u_i,v_i]\).  
Apply to the five pairs,
respectively, the following word substitutions
\[
\mathcal N_m=(\mathcal L_1,\ 
\mathcal U_{-1}\circ\mathcal L_{m-1},\ 
\mathcal L_{-1},\mathcal L_{-1},\mathcal U_{-1}).
\]
The five handle maps in $\mathcal N_m$ induce the matrices
$L_1,R_m,L_{-1},L_{-1},U_{-1}$ on the first homology, since the map
on the second pair is $\mathcal U_{-1}\circ\mathcal L_{m-1}$ and
\[
U_{-1}L_{m-1}
=\begin{pmatrix}1&-1\\0&1\end{pmatrix}
 \begin{pmatrix}1&0\\m-1&1\end{pmatrix}
=\begin{pmatrix}2-m&-1\\m-1&1\end{pmatrix}
=R_m.
\]

For adjacent pairs let \(C=[u_i,v_i]\) and define the based swap \(\mathcal H_i\) by
\[
\begin{aligned}
u_i&\longmapsto C u_{i+1}C^{-1},&
v_i&\longmapsto C v_{i+1}C^{-1},\\
u_{i+1}&\longmapsto u_i,&
v_{i+1}&\longmapsto v_i.
\end{aligned}
\]
It fixes all other generators and preserves
\([u_i,v_i][u_{i+1},v_{i+1}]\), since the new commutators are
\(C[u_{i+1},v_{i+1}]C^{-1}\) and \(C\). It is an invertible substitution of the free group preserving the full oriented boundary word $\prod_j[u_j,v_j]$. 
The Dehn-Nielsen-Baer theorem therefore realizes it by an orientation-preserving diffeomorphism of the genus-five surface with one boundary component fixed pointwise 
\cite[Theorem~3.5(1)]{Sakasai}. 
Capping that boundary gives a based diffeomorphism of the closed surface.  
Apply the adjacent swaps in the
chronological order \(\mathcal H_3,\mathcal H_2,\mathcal H_1,
\mathcal H_0,\mathcal H_1,\mathcal H_2\), and call the resulting
automorphism \(\mathcal H_\pi\). It sends source handle \(j\) to
destination handle \(\pi(j)\), with \(\pi=(3,1,2,4,0)\).

By the $(2,4)$ column of the first table in Subsection~\ref{subsec:matrices},
$\mathcal H_\pi\circ\mathcal N_m$ induces $S_m$ from the
two-sheeted multisection $F_{3,2}$ to the four-sheeted multisection $F_{\beta_k}$.
Our gluing $\psi_m:F_{\beta_k}\to F_{3,2}$ goes in the reverse
direction. In the above ordered surface generators, its word map is
$\mathcal N_m^{-1}\circ\mathcal H_\pi^{-1}$, and its matrix on
the first homology is $S_m^{-1}$. 
For the inverse swap, if
\(C'=[u_{i+1},v_{i+1}]\), then \(\mathcal H_i^{-1}\) sends
\((u_i,v_i)\) to \((u_{i+1},v_{i+1})\) and sends
\((u_{i+1},v_{i+1})\) to
\((C'^{-1}u_iC',C'^{-1}v_iC')\). 
Since $\mathcal U_c^{-1}=\mathcal U_{-c}$ and
$\mathcal L_c^{-1}=\mathcal L_{-c}$, the inverse substitutions
on the five handles, in the same order, are
\[
\mathcal N_m^{-1}
=\bigl(
  \mathcal L_{-1},\,
  \mathcal L_{1-m}\circ\mathcal U_1,\,
  \mathcal L_1,\,
  \mathcal L_1,\,
  \mathcal U_1
\bigr).
\]
For each fixed $m$, we use these same ten words for every $k$.

We may choose $\psi_m$ to preserve area while inducing the prescribed
map on the based fundamental groups. Indeed, let $\omega_1$ and
$\omega_2$ be the area forms on $F_{\beta_k}$ and $F_{3,2}$,
respectively, rescaled to have equal total area as in the proof of
Proposition~\ref{prop:four-sum}. Then $\psi_m^*\omega_2-\omega_1$
is exact. Choose a primitive $\lambda$ that vanishes at the source
basepoint $\widetilde y_0$. To do so, start with any primitive and
subtract $df$, where $f$ is a smooth function whose differential at
$\widetilde y_0$ agrees with the value of that primitive there.
For $\omega_t=(1-t)\omega_1+t\psi_m^*\omega_2$, the Moser vector
field defined by $\iota_{X_t}\omega_t=-\lambda$ therefore vanishes
at $\widetilde y_0$. Its flow $h_t$ fixes $\widetilde y_0$ and
satisfies $h_1^*\psi_m^*\omega_2=\omega_1$. Thus $\psi_m\circ h_1$
is area preserving and, since $h_t$ fixes the basepoint throughout,
induces the same ten words as $\psi_m$.

The $(2,4)$ pair calculation in Subsection~\ref{subsec:matrices} gives
\[
H_1(\pi_1(Z_{k,\psi_m});\mathbb Z)
 \cong (\mathbb Z/2)^2\oplus\mathbb Z/(4m).
\]

\begin{proposition}\label{prop:s3-distinction}
For the based gluings as above, we have
\[
\#\operatorname{Epi}(\pi_1(Z_{k,\psi_m}),S_3)
=\begin{cases}
36,&3\nmid k,\\
24,&3\mid k,
\end{cases}
\]
where $\#\operatorname{Epi}(G,H)$ denotes the cardinality of the set of surjective homomorphisms from $G$ to $H$.  
Here distinct homomorphisms are counted separately, even when they differ by
an automorphism of $S_3$.
\end{proposition}

\begin{proof}
Let $\rho:\pi_1(Z_{k,\psi_m})\to S_3$ be a homomorphism. We use the
same letters for the generators and their images under $\rho$.
The Artin relations imply that $x_1,x_2,x_3,x_4$ are conjugate,
and the meridian relation gives $t=x_1^{-1}$. We count the
epimorphisms according to whether $x_1$ is the identity, a
three-cycle, or a transposition.

We first determine which generators can have odd image. The
abelianization computed above is $(\Z/2)^2\oplus\Z/(4m)$,
so there are exactly eight homomorphisms to $\Z/2$. To obtain all
eight, assign arbitrary values to $a,\alpha,x_1$, give
$x_2,x_3,x_4,t$ the same value as $x_1$, and set
\begin{equation}\label{eq:8generators}
 b=c=d=\eta=\gamma=\delta=\epsilon=\zeta=0.
\end{equation}
These assignments satisfy the two complement presentations of
Example~\ref{ex:beta_k}. They also satisfy the gluing relations of
Proposition~\ref{prop:explicit-closed-group}. Indeed, $U_k$ has
total exponent sum $4k+4$, so $A=U_k^{-1}a^4$ has value zero in
$\Z/2$. Likewise $A'=\alpha^2$ has value zero. Thus every entry
of both surface tuples has value zero. Applying this observation
to the composite of $\rho$ with the sign map $S_3\to S_3/A_3\cong\Z/2$, we
see that $A,A'$ and the eight generators in (\ref{eq:8generators})
all lie in the alternating subgroup $A_3=\{1,(123),(132)\}$.

We claim that if $\rho$ is surjective, all ten of these elements
are the identity. To prove this, write $A_3\cong\Z/3$ additively,
and let $\sigma,\tau\in\{1,-1\}$ be the signs of $a,\alpha$,
respectively. Conjugation by an even permutation fixes each
element of $A_3$, and conjugation by a transposition takes each
element to its inverse. Thus conjugation by $a$ and $\alpha$
acts on $A_3$ by multiplication by $\sigma$ and $\tau$.

We now write the gluing relations in $A_3\cong\mathbb Z/3$.
Recall that an inverse swap exchanges two adjacent pairs and
conjugates the entries of one pair by the commutator of the
other. Since all entries lie in $A_3$, this commutator is
trivial. Thus, after applying the inverse handle substitutions
to the five pairs of $\mathcal L'$, we need only reorder the
resulting pairs. The forward gluing sends handle $i$ on the
two-sheeted side to handle $\pi(i)$ on the four-sheeted side.
Hence the pair in position $j$ of $\mathcal L_k$ is identified
with the transformed pair from position $\pi^{-1}(j)$ of
$\mathcal L'$. For $j=0,1,2,3,4$, these positions are
$4,1,2,0,3$, respectively.

For example, the first pair $(A,b)$ of $\mathcal L_k$ corresponds
to the pair $(E'_1,H'_1)$ in position $4$ of $\mathcal L'$.
Since conjugation by $\alpha$ acts on $A_3$ by multiplication
by $\tau$, this pair has image $(\tau\epsilon,\tau\zeta)$.
The inverse substitution on handle $4$ is $\mathcal U_1$, which
adds the first entry to the second. This gives the first row in (\ref{eq:relations})
below. For the second row, the pair in position $1$ is
$(\gamma,\delta)$, and its inverse substitution sends $(u,v)$
to $(uv^{1-m},vuv^{1-m})$. In additive notation its image is
therefore $(\gamma+(1-m)\delta,\gamma+(2-m)\delta)$.

On the four-sheeted side, $(C_r,D_r)=(a^rca^{-r},a^rda^{-r})$
has image $(\sigma^r c,\sigma^r d)$. Reading the transformed
pairs of $\mathcal L'$ in the order $4,1,2,0,3$ now gives all
five pairs of gluing relations:
\begin{equation}\label{eq:relations}
\begin{aligned}
(A,b)              &= \bigl(\tau\epsilon,\tau(\zeta+\epsilon)\bigr),\\
(c,d)              &= \bigl(\gamma+(1-m)\delta,\gamma+(2-m)\delta\bigr),\\
(\sigma c,\sigma d)&= (\epsilon+\zeta,\zeta),\\
(c,d)              &= (A'-\eta,\eta),\\
(\sigma c,\sigma d)&= \bigl(\tau(\gamma+\delta),\tau\delta\bigr).
\end{aligned}
\end{equation}
Subtracting the two entries in the second row gives $\delta=d-c$.
The third row gives $\zeta=\sigma d$ and $\epsilon=\sigma(c-d)=-\sigma\delta$, and the
last row gives $\delta=\sigma\tau d$. The first and fourth rows
therefore give
\[
 A=\tau\epsilon=-d,\qquad \eta=d.
\]

Suppose $d\ne0$. 
Since $d\in A_3$, \(d\) is one of the two three-cycles and \(A_3=\{1,d,d^{-1}\}\).
Since $\gamma,\delta,\epsilon,\zeta$ lie in the
abelian group $A_3$ by (\ref{eq:8generators}), the relation
$[\alpha,\eta][\gamma,\delta][\epsilon,\zeta]=1$ gives
$[\alpha,\eta]=1$. Thus $\alpha$ commutes with $d$, as does each
$x_j$ by the four-strand complement relations in (\ref{eq:4-strand-relations}). The centralizer
of a three-cycle in $S_3$ is $A_3$, so $\alpha$ and all four
$x_j$ lie in $A_3$. In particular $\alpha$ commutes with
$t=x_1^{-1}$. Since $\alpha t\alpha^{-1}=t^{-1}$, we have $t=1$,
and hence all four $x_j$ are the identity. It follows that $U_k=1$
and $a^4=A$ is a nontrivial three-cycle. This forces $a$ to be
even, since every odd permutation in $S_3$ has order two.
Thus $\rho$ sends every generator of $\pi_1(Z_{k,\psi_m})$ into
$A_3$. Therefore, if $\rho$ is surjective, we must have $d=0$.
The above equations now give
$\delta=c=\epsilon=\zeta=\eta=0$ and then $\gamma=b=A=A'=0$.
This proves the claim. We now return to multiplicative notation.

\smallskip
\noindent\emph{Case\/ $1$. The meridian $x_1$ is the identity.}
All four $x_j$ and $t$ are then the identity. The only generators
that can have nontrivial image are $a$ and $\alpha$, and the
remaining relations are $a^4=\alpha^2=1$. Each of $a,\alpha$
can therefore be the identity or a transposition. They generate
$S_3$ exactly when they are distinct transpositions. This gives
$3\cdot2=6$ epimorphisms.

\smallskip
\noindent\emph{Case\/ $2$. The meridian $x_1$ is a three-cycle.}
The $x_j$ lie in $A_3$ and hence commute. The Artin relations in (\ref{eq:4-strand-relations}) and (\ref{eq:Artin})
therefore reduce to $ax_ja^{-1}=x_{j+1}$, with indices taken
modulo four. If $a$ were even, all four $x_j$ would equal $x_1$,
contradicting $x_1x_2x_3x_4=1$. Hence $a$ is a transposition and
\begin{equation}\label{eq:x-tuple}
 (x_1,x_2,x_3,x_4)=(x_1,x_1^{-1},x_1,x_1^{-1}).
\end{equation}
Also $\alpha$ is a transposition, since it inverts the nontrivial
three-cycle $t$. (Note that $\alpha=a_1$ in (\ref{eq:2-strand-group}).) Conversely, choose either three-cycle for $x_1$
and any transpositions for $a$ and $\alpha$. The displayed tuple (\ref{eq:x-tuple})
satisfies the sphere and Artin relations. To check $A=1$, note
that $\Delta=x_1x_2=1$ and $z=x_1$, so the four factors of $U_k$
have images $x_1^{-1},x_1,x_1^{-1},x_1$, in that order.
Thus $U_k=1=a^4$. 
Taking the eight remaining base generators to
be the identity satisfies all remaining relations. Each map is
surjective since $x_1$ and $a$ generate $S_3$. This gives
$2\cdot3\cdot3=18$ epimorphisms.

\smallskip
\noindent\emph{Case\/ $3$. The meridian $x_1$ is a transposition.}
Fix one of the three choices for $x_1$. Since $t=x_1=t^{-1}$,
the relation $\alpha t\alpha^{-1}=t^{-1}$ implies
$\alpha x_1\alpha^{-1}=x_1$, so 
$\alpha\in\{1,x_1\}$. It remains to determine $a$ and the other
three meridians. 
Label the three letters by the elements of $\mathbb F_3=\{0,1,2\}$, taking
the point fixed by $x_1$ to be $0$. The six maps
$y\mapsto u+\sigma y$, with $u\in\mathbb F_3$ and $\sigma=\pm1$,
are distinct permutations and hence exhaust $S_3$. 
Each map $y\mapsto u-y$ fixes the unique solution of $2y=u$ in
$\mathbb F_3$ and interchanges the other two points, so the three
choices of $u$ give the three transpositions.
Thus we may write
\[
 x_j(y)=r_j-y,\qquad a(y)=u+\sigma y,
\]
where $r_j\in\mathbb F_3$ and $r_1=0$. 
Set $r=r_2$.
Conjugation by $a$ sends the reflection $y\mapsto v-y$ to the
reflection with parameter $2u+\sigma v$. The sphere relation gives
$r_4=r_3-r$. Also $\Delta=x_1x_2$ is translation by $-r$, so
$z=\Delta^kx_1\Delta^{-k}$ is the reflection with parameter
$-2kr=kr$ in $\mathbb F_3$.

We now use the Artin relations $a x_j a^{-1}=\Phi(x_j)$ for
$j=1,4,2$ in Example~\ref{ex:beta_k}. 
On the left, conjugation by $a$ sends the reflection
$x_j(y)=r_j-y$ to the reflection with parameter $2u+\sigma r_j$.
On the right, $\Phi(x_1)=\Delta z^{-1}$ has parameter $(k-1)r$,
since $z^{-1}=z$ and $\Delta$ subtracts $r$. Also,
$\Phi(x_4)=z$ has parameter $kr$. Finally,
$\Phi(x_2)=z x_3 z^{-1}$ has parameter $2kr-r_3$, by the same
conjugation formula applied to $z(y)=kr-y$. Equating the
parameters on the two sides and using $r_1=0$ and $r_2=r$, we obtain
\[
2u=(k-1)r,\qquad
2u+\sigma r_4=kr,\qquad
2u+\sigma r=2kr-r_3.
\]

The first two give $r_4=\sigma r$ and $r_3=(1+\sigma)r$.
Substitution into the third gives 
\begin{equation}\label{eq:k-2sigma}
(k-2\sigma)r=0. 
\end{equation}
The Artin relation at $x_3$ then follows from the other three because
$\Phi_{\beta_k}$ preserves $x_1x_2x_3x_4$. These conditions
also suffice for the gluing. 
By the Artin relations and the definition of $U_k$, we have
\[
a^4x_1a^{-4}=\Phi^4(x_1)=U_kx_1U_k^{-1}.
\]
It follows that $A=U_k^{-1}a^4$ centralizes $x_1$. Since each
$x_j$ is a transposition and the exponent sum of $U_k$ is $4k+4$,
the permutation $U_k$ is even. The fourth power $a^4$ is also
even, so $A$ is even. Any permutation commuting with $x_1$ must
fix its unique fixed letter, so the centralizer of $x_1$ in $S_3$
is $\{1,x_1\}$. Hence $A=1$. With $b=c=d=1$, every entry of
$\mathcal L_k$ is therefore trivial. On the two-sheeted side,
$\alpha\in\{1,x_1\}$ gives $A'=\alpha^2=1$. With
$\eta=\gamma=\delta=\epsilon=\zeta=1$, every entry of
$\mathcal L'$ is trivial as well. Thus all ten surface gluing
relations hold.

If $r=0$, then $u=r_3=r_4=0$, and every generator lies in
$\{1,x_1\}$. If $r\ne0$, the distinct transpositions $x_1$ and
$x_2$ generate $S_3$. In this case $k=2\sigma$ in $\mathbb F_3$ 
by (\ref{eq:k-2sigma}).
There is no such $\sigma$ when $3\mid k$, and there is exactly one
$\sigma$ when $3\nmid k$. For that $\sigma$, either nonzero value of $r$
determines $u,r_3,r_4$ uniquely. Including the three choices for
$x_1$ and the two for $\alpha$, we obtain $3\cdot2\cdot2=12$
epimorphisms when $3\nmid k$ and none when $3\mid k$.

Adding the three cases gives $6+18+12=36$ when $3\nmid k$ and
$6+18=24$ when $3\mid k$. 
\end{proof}

\begin{proof}[Proof of Theorem~{\rm \ref{thm:main-fixed-gluing}}]
Take $k=1,3$ in the gluings constructed in this section. Proposition~\ref{prop:four-sum} gives the symplectic sums and their Euler characteristic and signature. Proposition~\ref{prop:aspherical-sums} gives asphericity and minimality. The homology calculation in Subsection~\ref{subsec:matrices} gives their common first homology, and Proposition~\ref{prop:s3-distinction} distinguishes their fundamental groups. The rational cohomology and eventual noncomplexity follow from Theorem~\ref{thm:noncomplex}, applied to the two sequences.
\end{proof}

\section{An infinite family of genus-two fibered knots}

We now apply the hyperelliptic quotient construction of
Section~\ref{sec:genus-two} to Misev's genus-two fibered knots
$K_n$, $n\geq0$ in \cite{Misev}. These knots have the same Seifert form as
$K_0=T(2,5)$. For $n>0$, their knot-fiber monodromies are
pseudo-Anosov and have pairwise distinct stretch factors.
We use the common Seifert form to determine the degrees of the
two components of the branch locus.

Let $\widehat\varphi_n$ denote the monodromy on the closed surface
obtained by capping the fiber of $K_n$. Its mapping torus is the
zero surgery $M_{K_n}$. 
By the genus-two case of the Birman--Hilden theorem~\cite{BirmanHilden},
the capped monodromy determines a mapping class of the sphere with
six marked branch points. As in Section~\ref{sec:four-strand}, choose
a representative $\bar\varphi_n$ that is a rotation in local
coordinates near each branch point. Of its two smooth lifts to the
genus-two surface, take the one isotopic to $\widehat\varphi_n$ as
the representative of the capped monodromy.

The fiberwise quotient of $M_{K_n}$
is the mapping torus of $\bar\varphi_n$. Since $\bar\varphi_n$ preserves
orientation, its mapping torus is diffeomorphic to $S^1\times S^2$.
Taking the product with another circle gives a double branched cover
\[
M_{K_n}\times S^1\longrightarrow T^2\times S^2.
\]

\begin{proposition}\label{prop:five}
For each $n\geq0$, the branch locus consists of a section $B_{1,n}$
and a connected degree-five multisection $B_{5,n}$. The latter is a
torus of self-intersection zero representing the homology class\/ 
$5[T^2\times\{\mathrm{pt}\}]$.
\end{proposition}

\begin{proof}
The characteristic polynomial of the action of $\widehat\varphi_n$
on integral first homology is the Alexander polynomial of $K_n$,
normalized to be monic with constant term one. Since the Seifert
form agrees with that of $T(2,5)$, this polynomial is
$t^4-t^3+t^2-t+1$. 
To determine how $\bar\varphi_n$ permutes the six branch points,
we use their equivariant correspondence with the six odd quadratic
refinements on $H_1(\Sigma_2;\mathbb F_2)$, as in
Section~\ref{sec:genus-two}. 

Let $M$ denote the induced action of $\widehat\varphi_n$ on
$H_1(\Sigma_2;\mathbb F_2)$. Its characteristic polynomial satisfies
\[
(t+1)(t^4-t^3+t^2-t+1)=t^5+1=t^5-1
\quad\text{in }\mathbb F_2[t].
\]
By Cayley-Hamilton theorem, $M^5=I$. 
The characteristic polynomial of $I$ is $(t-1)^4=t^4+1$, so $M$ has order five. As in the proof of Lemma~\ref{lem:three-three} in 
Section~\ref{sec:genus-two}, the six branch points correspond to
the six odd quadratic refinements. The action on these refinements
identifies $\operatorname{Sp}(4,\mathbb F_2)$ with $S_6$. Thus
$\bar\varphi_n$ acts on the branch points by a five-cycle with one
fixed point.

The fixed point traces a circle in the sphere mapping torus, and
the other five points trace a single circle covering the base
circle five times. Crossing these circles with the second $S^1$
gives the disjoint tori $B_{1,n}$ and $B_{5,n}$. For the product
identification, choose an isotopy from $\bar\varphi_n$ to the
identity fixing its fixed branch point. Such an isotopy exists
since the mapping class group of the sphere with one marked point
is trivial. Then $B_{1,n}$ is a constant section. Write
$A=[T^2\times\{\mathrm{pt}\}]$ and
$F=[\{\mathrm{pt}\}\times S^2]$. Since $B_{5,n}$ has degree five,
$[B_{5,n}]=5A+dF$ for some $d\in\Z$. Intersecting with $B_{1,n}$
gives $d=0$, and $A^2=0$ gives $[B_{5,n}]^2=0$.
\end{proof}

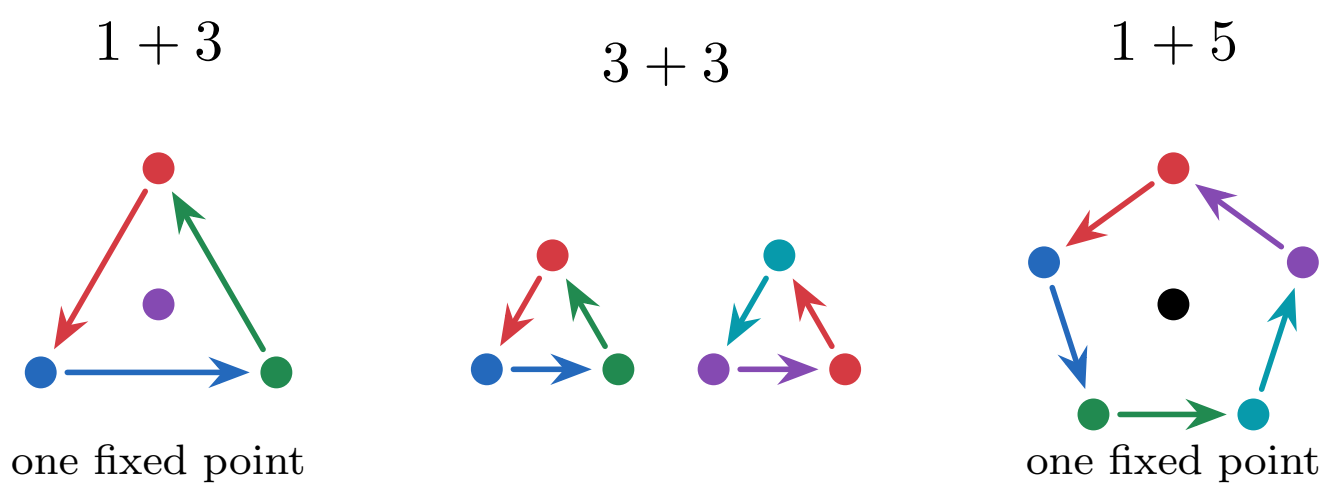
\begin{figure}[htbp]
\centering
\begin{tikzpicture}[x=.75cm,y=.75cm,>=Stealth,line cap=round]
\node at (0,2.35) {$1+3$};
\draw[->,color=strandred,line width=1.05pt,shorten >=5pt,shorten <=5pt] (0.0000,1.2500) -- (-1.0392,-0.5500);
\draw[->,color=strandblue,line width=1.05pt,shorten >=5pt,shorten <=5pt] (-1.0392,-0.5500) -- (1.0392,-0.5500);
\draw[->,color=strandgreen,line width=1.05pt,shorten >=5pt,shorten <=5pt] (1.0392,-0.5500) -- (0.0000,1.2500);
\fill[strandred] (0.0000,1.2500) circle (3pt);
\fill[strandblue] (-1.0392,-0.5500) circle (3pt);
\fill[strandgreen] (1.0392,-0.5500) circle (3pt);
\fill[twistpurple] (0,.05) circle (3pt);
\node[font=\scriptsize] at (0,-1.35) {one fixed point};
\path (-1.9,-1.5) rectangle (1.9,2.5);
\end{tikzpicture}
\hspace{.25cm}
\begin{tikzpicture}[x=.75cm,y=.75cm,>=Stealth,line cap=round]
\node at (0,2.35) {$3+3$};
\draw[->,color=strandred,line width=1.05pt,shorten >=5pt,shorten <=5pt] (-1.0000,0.6700) -- (-1.5802,-0.3350);
\draw[->,color=strandblue,line width=1.05pt,shorten >=5pt,shorten <=5pt] (-1.5802,-0.3350) -- (-0.4198,-0.3350);
\draw[->,color=strandgreen,line width=1.05pt,shorten >=5pt,shorten <=5pt] (-0.4198,-0.3350) -- (-1.0000,0.6700);
\fill[strandred] (-1.0000,0.6700) circle (3pt);
\fill[strandblue] (-1.5802,-0.3350) circle (3pt);
\fill[strandgreen] (-0.4198,-0.3350) circle (3pt);
\draw[->,color=strandturquoise,line width=1.05pt,shorten >=5pt,shorten <=5pt] (1.0000,0.6700) -- (0.4198,-0.3350);
\draw[->,color=twistpurple,line width=1.05pt,shorten >=5pt,shorten <=5pt] (0.4198,-0.3350) -- (1.5802,-0.3350);
\draw[->,color=strandred,line width=1.05pt,shorten >=5pt,shorten <=5pt] (1.5802,-0.3350) -- (1.0000,0.6700);
\fill[strandturquoise] (1.0000,0.6700) circle (3pt);
\fill[twistpurple] (0.4198,-0.3350) circle (3pt);
\fill[strandred] (1.5802,-0.3350) circle (3pt);
\path (-1.9,-1.5) rectangle (1.9,2.5);
\end{tikzpicture}
\hspace{.25cm}
\begin{tikzpicture}[x=.75cm,y=.75cm,>=Stealth,line cap=round]
\node at (0,2.35) {$1+5$};
\draw[->,color=strandred,line width=1.05pt,shorten >=5pt,shorten <=5pt] (0.0000,1.2500) -- (-1.1413,0.4208);
\draw[->,color=strandblue,line width=1.05pt,shorten >=5pt,shorten <=5pt] (-1.1413,0.4208) -- (-0.7053,-0.9208);
\draw[->,color=strandgreen,line width=1.05pt,shorten >=5pt,shorten <=5pt] (-0.7053,-0.9208) -- (0.7053,-0.9208);
\draw[->,color=strandturquoise,line width=1.05pt,shorten >=5pt,shorten <=5pt] (0.7053,-0.9208) -- (1.1413,0.4208);
\draw[->,color=twistpurple,line width=1.05pt,shorten >=5pt,shorten <=5pt] (1.1413,0.4208) -- (0.0000,1.2500);
\fill[strandred] (0.0000,1.2500) circle (3pt);
\fill[strandblue] (-1.1413,0.4208) circle (3pt);
\fill[strandgreen] (-0.7053,-0.9208) circle (3pt);
\fill[strandturquoise] (0.7053,-0.9208) circle (3pt);
\fill[twistpurple] (1.1413,0.4208) circle (3pt);
\fill[black] (0,.05) circle (3pt);
\node[font=\scriptsize] at (0,-1.35) {one fixed point};
\path (-1.9,-1.5) rectangle (1.9,2.5);
\end{tikzpicture}
\caption{The cycles of the branch-point permutation for the genus-one
knots, the even Stallings twists of the square knot, and Misev's
genus-two knots, respectively.}
\label{fig:branch-orbits}
\end{figure}

Figure~\ref{fig:branch-orbits} shows the branch-point orbits in
the three knot constructions considered in this paper.

Write $T^2=S^1\times S^1$, where the first factor is the base
circle of the sphere mapping torus. Choose a smooth map
$u:\Sigma_2\to S^1$ such that
$u_*:\pi_1(\Sigma_2)\to\pi_1(S^1)\cong\Z$ is surjective.
The map $x\mapsto (u(x),1)$ pulls the product sphere bundle
$T^2\times S^2$ back to $\Sigma_2\times S^2$. In this pullback,
the multisection $B_{5,n}$ gives the surface
\[
F_{5,n}
=\bigl\{(x,z)\in\Sigma_2\times S^2
  \mid ((u(x),1),z)\in B_{5,n}\bigr\}.
\]
Thus the five points of $F_{5,n}$ over $x$ are precisely the
five points of $B_{5,n}$ over $(u(x),1)$. In particular, the
projection $F_{5,n}\to\Sigma_2$ is an unbranched covering of
degree five.

To see that $F_{5,n}$ is connected, choose a loop $\gamma$ based
at $y_0\in\Sigma_2$ such that $u\circ\gamma$ represents the
positive generator of $\pi_1(S^1,u(y_0))$. Such a loop exists
because $u_*$ is surjective. Following the five points of
$F_{5,n}$ around $\gamma$ gives the same permutation as following
the five points of $B_{5,n}$ once around the first circle of
$T^2$, with the second coordinate fixed at $1$. By
Proposition~\ref{prop:five}, this permutation is a five-cycle.
Hence the lifts of successive powers of $\gamma$ starting at any
one point over $y_0$ reach all five points over $y_0$. The covering
is therefore connected.

$F_{5,n}$ is disjoint from the pulled-back
constant section. As in the proof of Proposition~\ref{prop:five},
its degree and its intersection with that section give
\[
[F_{5,n}]=5[\Sigma_2\times\{\mathrm{pt}\}],\qquad [F_{5,n}]^2=0.
\]
Its Euler characteristic is $5\chi(\Sigma_2)=-10$, so it has
genus six. Since the projection restricts to an orientation-preserving
covering on $F_{5,n}$, adding a sufficiently large multiple of the
base area form to the fiber area form makes $F_{5,n}$ symplectic.

\begin{proposition}\label{prop:twelve}
For every $n\geq0$, two copies of $\Sigma_2\times S^2$ admit a
minimal symplectic sum along $F_{5,n}$ and $F_{5,0}$ with
\[
e=12,\qquad \sigma=0,\qquad b_1=0.
\]
\end{proposition}

\begin{proof}
Let $X=\Sigma_2\times S^2$. Both multisections have genus six and
trivial normal bundle. Choose the product symplectic forms so that
the two surfaces have the same area. For any symplectic gluing,
additivity of the signature gives $\sigma=0$, and
\[
e=2e(X)-2\chi(\Sigma_6)=2(-4)-2(-10)=12.
\]
It only remains to choose the surface identification so that the first Betti number vanishes.

Write $F_1=F_{5,n}$ and $F_2=F_{5,0}$ in the respective copies of
$X$, and let $q_i:F_i\to\Sigma_2$ be the covering projections.
Let $M_i=X\setminus\operatorname{Int}\nu F_i$ be the exterior of
$F_i$. By excision and the Thom isomorphism for the oriented normal
bundle of $F_i$, we have
\[
H_k(X,M_i;\Q)\cong H_{k-2}(F_i;\Q).
\]
Since $F_i$ is connected, $H_2(X,M_i;\Q)\cong\Q$, with generator
represented by an oriented normal disk to $F_i$. Moreover,
$H_1(X,M_i;\Q)=0$. Thus the relevant part of the exact sequence
of the pair $(X,M_i)$ is
\[
H_2(X;\Q)\longrightarrow \Q
\xrightarrow{\;\partial\;} H_1(M_i;\Q)
\longrightarrow H_1(X;\Q)\longrightarrow 0.
\]
The first map sends a class to its algebraic intersection number
with $F_i$. Since $q_i:F_i\to\Sigma_2$ is an orientation-preserving
covering of degree five, a sphere fiber meets $F_i$ transversely
in five points, all with intersection sign $+1$. Thus the fiber
class maps to $5$, and the first map is onto over $\Q$. By
exactness, $\partial=0$. The boundary of the normal disk is a
meridian of $F_i$, so the meridian represents zero in
$H_1(M_i;\Q)$. Exactness also shows that inclusion induces an
isomorphism $H_1(M_i;\Q)\to H_1(X;\Q)$. Composing with the
projection $X=\Sigma_2\times S^2\to\Sigma_2$ identifies
$H_1(M_i;\Q)$ with $H_1(\Sigma_2;\Q)$.

Choose a normal framing and write $\partial M_i=F_i\times S^1$,
where the circle factor is a meridian. Inside $\nu F_i$, the
inclusion of $\partial M_i$ is homotopic to projection onto $F_i$
followed by inclusion of the zero section. Projecting the zero
section to $\Sigma_2$ gives $q_i$. Hence, under the identification
above, the boundary inclusion induces the map
\[
H_1(F_i;\Q)\oplus\Q\longrightarrow H_1(\Sigma_2;\Q),
\qquad (a,r)\longmapsto(q_i)_*a.
\]

Let $Z$ be the sum obtained by a boundary gluing covering
$\psi:F_1\to F_2$. In the Mayer-Vietoris sequence for
$Z=M_1\cup M_2$, the map from the first homology of the common
boundary is the difference of the two inclusion maps. A class
with surface component $a$ on the first boundary has surface
component $\psi_*a$ on the second. Any meridian component maps
to zero in either exterior. The image of the Mayer-Vietoris map
is therefore the image of
\[
\big((q_1)_*,-(q_2)_*\psi_*\big):
H_1(F_1;\Q)\longrightarrow
H_1(\Sigma_2;\Q)\oplus H_1(\Sigma_2;\Q).
\]
Finally, the common boundary and both exteriors are connected, so
the map on $H_0$ sends $1$ to $(1,-1)$ and is injective.
Exactness now identifies $H_1(Z;\Q)$ with the cokernel of the
displayed map.

To make this map surjective, it suffices to make its dual injective.
Put $W_i=q_i^*H^1(\Sigma_2;\Q)\subset H^1(F_i;\Q)$. Since $q_i$
has degree five, the cup product pairing on pullbacks is five times
the pairing on the base. Thus $q_i^*$ is injective and $W_i$ is a
four-dimensional symplectic subspace. The dual map sends $(a,b)$ to
$q_1^*a-\psi^*q_2^*b$, so it is injective precisely when
\[
W_1\cap\psi^*W_2=0.
\]
The space $W_1^{\perp}$ has dimension eight
and contains a four-dimensional symplectic subspace. 
Choose a symplectic basis of $W_2$ and a symplectic basis of a
four-dimensional symplectic subspace of $W_1^\perp$. Extending
these to symplectic bases of $H^1(F_2;\Q)$ and $H^1(F_1;\Q)$,
respectively, gives a rational symplectic isomorphism
\[
R:H^1(F_2;\Q)\longrightarrow H^1(F_1;\Q)
\]
that sends the first basis to the second. Thus
$R(W_2)\subset W_1^\perp$, and hence $W_1\cap R(W_2)=0$.

Choose integral symplectic bases for $H^1(F_1;\Z)$ and
$H^1(F_2;\Z)$, and fix an integral basis for $H^1(\Sigma_2;\Z)$.
For $i=1,2$, let $P_i$ be the $12\times4$ matrix of $q_i^*$
in these bases. Its four columns form a basis for $W_i$ over $\Q$.
For any rational symplectic isomorphism
$S:H^1(F_2;\Q)\to H^1(F_1;\Q)$, write $S$ also for its matrix
in the chosen bases. Then the first four columns of
\[
\bigl[\,P_1\mid S P_2\,\bigr]
\]
form a basis for $W_1$, and the last four form a basis for
$S(W_2)$. Thus the matrix has rank eight precisely when
$W_1\cap S(W_2)=0$. Taking $S=R$, the isomorphism constructed
above, we have $W_1\cap R(W_2)=0$. We may therefore choose eight
rows of $\bigl[\,P_1\mid R P_2\,\bigr]$ for which the resulting
$8\times8$ matrix has nonzero determinant. Let $D(S)$ denote
the determinant obtained by taking those same eight rows of
$\bigl[\,P_1\mid S P_2\,\bigr]$. Then $D(R)\ne0$, and it
suffices to find an \emph{integral}\/ symplectic matrix $S$ for which
$D(S)\ne0$.

In the chosen bases, regard $R$ as an
element of $\operatorname{Sp}(12,\Q)$. As in the proof of
Proposition~\ref{prop:unbounded-gluing}, write
\[
R=\tau_{v_1}(r_1)\cdots\tau_{v_\ell}(r_\ell),
\qquad
\tau_v(r)(x)=x+r\langle x,v\rangle v,
\]
where $v_j\in\Z^{12}$ and $r_j\in\Q$. 
Here we may take each $v_j$
to be integral because replacing $v$ by $dv$ and $r$ by $r/d^2$
gives the same transvection. Choose a positive integer $N$ such
that $Nr_j\in\Z$ for every $j$, and set
\[
S(t)=\tau_{v_1}(Nr_1t)\cdots\tau_{v_\ell}(Nr_\ell t).
\]
Each factor has entries in $\Z[t]$ and is an integral symplectic
matrix whenever $t$ is an integer. Hence $S(m)$ belongs to
$\operatorname{Sp}(12,\Z)$ for every integer $m$. Moreover,
$S(1/N)=R$. It follows that $D(S(t))$ is a polynomial in $t$
whose value at $1/N$ is $D(R)\ne0$. A nonzero polynomial has
only finitely many roots, so we may choose an integer $m$ with
$D(S(m))\ne0$. For this choice, $W_1\cap S(m)W_2=0$.

By the surjectivity of the symplectic representation of the mapping
class group~\cite[Theorem~6.4]{FarbMargalit}, there is an
orientation-preserving diffeomorphism $\psi:F_1\to F_2$ whose
pullback $\psi^*:H^1(F_2;\Z)\to H^1(F_1;\Z)$ has matrix $S(m)$.
Let $\omega_i$ denote the area form on $F_i$. Since the two
surfaces have equal area, $\omega_1$ and $\psi^*\omega_2$ have
the same integral over $F_1$. By Moser's theorem, there is a
diffeomorphism $h:F_1\to F_1$, isotopic to the identity, such
that $h^*(\psi^*\omega_2)=\omega_1$. Replacing $\psi$ by
$\psi\circ h$ therefore gives an area-preserving gluing with the
same action on cohomology. Since $W_1\cap\psi^*W_2=0$, the dual
of the Mayer--Vietoris map is injective. The map on first homology
is consequently surjective, so the symplectic sum has $b_1=0$.
Minimality follows from Proposition~\ref{prop:aspherical-sums}.
\end{proof}

\subsection*{Acknowledgments.}
The second author was partially supported by an NSERC Discovery grant 
RGPIN--2024--04401.
The authors used OpenAI's ChatGPT to help revise the exposition,
clarify steps in several proofs, and prepare figures and
computational checks.

\end{document}